\documentclass[11pt,oneside]{article}
\usepackage{amsmath,amsfonts,amssymb,amsthm}
\usepackage[colorlinks=true,linkcolor=RawSienna,citecolor=RawSienna,urlcolor=RawSienna,pdfstartview=FitB]{hyperref}
\usepackage[left=1.1in,right=1.1in,top=1.1in,bottom=1.1in]{geometry}
\usepackage[onehalfspacing]{setspace}
\usepackage{natbib}
\usepackage{hyphenat}
\usepackage{footmisc}
\usepackage{graphicx,float}
\usepackage[dvipsnames]{xcolor}
\usepackage{appendix}
\newtheorem{theorem}{Theorem}[]

\newtheorem{lemma}{Lemma}[]

\newtheorem{proposition}{Proposition}[]

\numberwithin{equation}{section}
\numberwithin{figure}{section}
\numberwithin{table}{section}

\newcommand{\rbk}[1]{\left(#1\right)}

\newcommand{\ms}{\boldsymbol}
\newcommand{\mb}{\mathbf}

\newcommand{\mbb}{\mathbb}

\begin{document}

\title{\vspace{-1.2cm}A strong law of large numbers related to multiple testing
Normal means}
\author{Xiongzhi Chen\thanks{Corresponding author: Department of Mathematics and Statistics, Washington State University,
Pullman, WA 99164, USA; Email: \texttt{xiongzhi.chen@wsu.edu}.}
\ and Rebecca W. Doerge\thanks{Office of the Dean, Mellon College of Science, 4400 Fifth Avenue, Pittsburgh, PA 15213, USA; Email: \texttt{rwdoerge@andrew.cmu.edu}.} }
\date{}
\maketitle

\begin{abstract}
Assessing the stability of a multiple testing procedure under dependence is important but very challenging. Even for multiple testing which among a set of Normal random variables have mean zero, which we refer to as the ``Normal means problem'', to date there lacks a classification of the type of dependence under which the strong law of large numbers (SLLN) holds for the numbers of rejections and false rejections. We introduce the concept of ``principal correlation structure (PCS)" that characterizes the type of dependence for which such SLLN holds, and establish the law. Further, we show that PCS ensures the SLLN for the false discover proportion when there is always a positive proportion of zero Normal means.
We also investigate the stability of two conditional multiple testing procedures for the Normal means problem, and show that the associated SLLN holds when in addition the decomposition of the covariance matrix of the Normal random variables that induces PCS is homogeneous in certain sense.
Our results also provide a formal way to check if the ``weak dependence'' assumption, a widely used assumption in the multiple testing literature, holds for the Normal means problem.
As by-products, we establish a universal bound on Hermite polynomials and a universal comparison result on the covariance of the indicator functions of the two p-values of testing the marginal means of a bivariate Normal random vector and the correlation between the two components of the vector. These are of their own interests.

\end{abstract}

\section{Introduction}

\label{secIntro}

We consider the ``Normal means problem under dependence", where one assesses which among many dependent Normal random variables have zero means. This problem has been widely encountered in
gene expression studies \cite{Owen:2005}, genome wide association studies \cite{Fan:2012}, and brain imaging analysis
\cite{Schwartzman:2015}, where test statistics are dependent and considered to be Normally distributed and whether their means are zero is tested.
Due to dependence, the behavior of the number of
rejections and the false discovery proportion (FDP, \cite{Genovese:2002}) of
a marginal multiple testing procedure (MTP) based on marginal observations
is unstable and sometimes even unpredictable; see, e.g., \cite{Finner:2007}, \cite{Owen:2005}, and \cite{Schwartzman2011}.
On the other hand, it has been observed that often a major part of the dependence among the observations of a given data set is induced
by some factors or latent variables. So, to adjust for dependence a
conditional multiple testing approach based on approximate factor models has
been taken. Specifically,  via the spectral decomposition of the covariance or correlation matrix of the Normal test statistics, this approach decomposes complicated dependence into a major part that is commonly referred to as ``(principal) factors", and applies MTPs to p-values of the test statistics conditional on the factors; see \cite{desai2012cross}, \cite{Fan:2017}, \cite{Fan:2012}, \cite{Friguet:2009} and \cite{Leek:2008}. However, an accurate estimate of the covariance or correlation matrix is usually needed to implement the approach.

This motivates us to investigate the stability of the marginal MTP and the conditional MTP for the Normal means problem. In particular, we address the following three important questions: Under what type of dependence will the SLLN hold for the random processes associated with the marginal MTP, i.e., the number of rejections, that of false rejections and the FDP, so that conditional multiple testing is not needed? Under what type of dependence will the SLLN hold for the same processes associated with the conditional MTP mentioned above? Is there a universal quantity that characterizes such dependence?

\subsection{Main contributions}

In this article, we provide partial answers to these questions. We classify the type of dependence for which the SLLN holds for the random processes associated with the marginal and conditional MTPs, and rigorously establish the law under the classified type of dependence.
Specifically, we propose the concept of ``principal correlation structure (PCS)'' (see definition (\ref{eqPCSdef})) and show that it alone is sufficient to ensure the SLLN for the number of rejections and the number of false rejections of the marginal MTP. Further, we show that the FDP of the marginal MTP satisfies the SLLN under PCS when there is always a positive proportion of zero Normal means. Since PCS is defined via the $l_1$-norm $\Vert \mathbf{R} \Vert_{1}$ of the correlation matrix $\mathbf{R}$ of the Normal random variables, $\Vert \mathbf{R} \Vert_{1}$ is a universal quantity that classifies the type of dependence for which the SLLN holds for the marginal MTP.
Interpreted slightly differently, for the Normal means problem, our result provides a fast and rigorous check on when a conditional MTP is needed, and converts assessing if the SLLN holds for the marginal MTP into testing the order of $\Vert \mathbf{R} \Vert_{1}$. Further, since the Normal means problem only involves variances of Normal random variables, our finding eliminates the need to accurately estimate the covariance or correlation matrix of these variables whenever a conditional MTP is not needed.

On the other hand, the SLLN associated with the marginal MTP is related to the ``weak dependence" assumption proposed by \cite{Storey:2004} that has been widely used in the multiple testing literature. The assumption requires that the SLLN hold for the number of rejections and that of false rejections. Up till now, there does not seem to be a formal way to check if the assumption holds. In view of this, our work fills this gap, and can be quite useful in studying the asymptotic conservativeness of an MTP under weak dependence.

For the conditional MTP based on an additive decomposition of the Normal random vector (see \eqref{eqModel}) that embraces the conditional multiple testing approach mentioned earlier, we show that the number of rejections and the number of false rejections satisfy the SLLN under PCS when additionally the decomposition is homogeneous (in the sense defined by \autoref{ThmSLLNStrongest}). So, our result classifies the type of dependence under which the SLLN holds for the conditional MTP via a universal quantity, the $l_1$-norm of the correlation matrix of the minor vector in the decomposition (see (\ref{eqModel}) and (\ref{eqPCSdef})). It also reveals that to retain the accuracy and stability of conditional multiple testing studied and advocated by \cite{desai2012cross}, \cite{Fan:2012}, \cite{Friguet:2009} and \cite{Leek:2008}, the residual dependence among the tests after adjusting for the factors or latent variables cannot change much as the number of test varies. In particular, for conditional multiple testing of \cite{desai2012cross} and \cite{Fan:2012}, the approximate factor model induced by the spectral decomposition of the covariance matrix of the Normal random variables cannot change much as the matrix changes.

The key difficulty in establishing the SLLN lies in connecting the covariance or correlation structure of a Normal random vector to that of the p-values induced by their corresponding tests on the Normal means. To deal with this, we derive
a universal comparison result for the correlation of a bivariate Normal random vector and the covariance of the indicator functions of the p-values of testing the means of its two marginal distributions. This is obtained by Mehler expansion and establishing a universal upper bound on Hermite polynomials. The comparison result and upper bound are of their own interests.

Our numerical study suggests that, compared to the method of ``principal factor approximation (PFA)" of \cite{Fan:2012}, constructing approximate factor models according to PCS
can achieve a much more adaptive decomposition of the dependence structure among the Normal random variables and give more stable conditional inferential results.

\subsection{Related work}

There are two works that are closely related to our work on studying the limiting behavior of random processes associated with an MTP. Firstly, when components of the Normal random vector all have variance $1$, \cite{Delattre:2016} proved that upon proper scaling and centering, the number of rejections converges in distribution to a Gaussian process, and obtained the central limit theorem for the FDP of the marginal MTP under the two-group model where the non-zero Normal means are identical. This implies the weak law of large numbers for the FDP under these settings. We remark that their results were obtained under a different set of conditions on the correlation matrix of the Normal random vector than those will be used in this work.

Secondly, \cite{Schwartzman:2015} introduced the concept of ``weak correlation'' among components of a Normal random vector.
 The authors of \cite{Schwartzman:2015} showed that, under the global null and when components of the Normal random vector all have variance $1$, the empirical cumulative distribution function (CDF) of the rejections of the marginal MTP converges in $L_2$-norm to the CDF of the standard Normal random variable if and only if components of the vector are weakly correlated. It is easy to see that weak correlation is weaker than PCS to be introduced later.

\subsection{Organization of article}

The rest of the article is organized as follows. \autoref{secModelPCS} introduces
marginal and conditional multiple testing Normal means and the concept of PCS.
\autoref{secNormalMeans} presents our theoretical results for the Normal means
problem under dependence.
\autoref{SecExample} provides an example for which the SLLN fails without PCS. \autoref{secSim} provides a simulation study on our theoretical findings.
\autoref{SecDiscussion} concludes the article with a discussion.
We provide in the main text short proofs but relegate longer ones to the appendices.
\autoref{SecConnection} discusses the relationship between PCS and the PFA of \cite{Fan:2012}.

\section{Multiple testing Normal means}

\label{secModelPCS}

Let $\left(  \Omega,\mathcal{F},\mathbb{P}\right)  $ be the probability space
on which all random vectors are defined, where $\Omega$ is the sample space,
$\mathcal{F}$ a sigma-algebra on $\Omega$, and $\mathbb{P}$ the probability
measure on $\mathcal{F}$. We will introduce marginal testing in \autoref{secMarginalTesting}, conditional testing
in \autoref{SecCondMTP}, and the concept of PCS in \autoref{SecDefPCS}.

Let $\mathsf{N}_{m}\left(  \mathbf{a},\mathbf{C}\right)  $ denote the Normal
distribution (and its density) with mean $\mathbf{a}\in\mathbb{R}^{m}$ and
covariance matrix $\mathbf{C}$. Pick $\boldsymbol{\mu}=\left(  \mu_{1},\ldots,\mu_{m}\right)  ^{\top}$
and set $\boldsymbol{\zeta}=\left(  \zeta_{1},\ldots,\zeta_{m}\right)
^{\top}$, $\boldsymbol{\eta}=\left(\eta_1,\ldots,\eta_m\right)^{\top}$ and $\mathbf{v}=\left(v_1,\ldots,v_m\right)^{\top}$. Suppose $\boldsymbol{\eta}\sim\mathsf{N}%
_{m}\left(  \mathbf{0},\boldsymbol{\Sigma}_{\boldsymbol{\eta}}\right)  $ and
$\mathbf{v}\sim\mathsf{N}_{m}\left(  \mathbf{0},\boldsymbol{\Sigma
}_{\mathbf{v}}\right)  $ and that $\boldsymbol{\eta}$ and $\mathbf{v}$ are
independent. Then setting $\boldsymbol{\zeta}=\boldsymbol{\mu
}+\boldsymbol{\eta}+\mathbf{v}$ gives $\boldsymbol{\zeta}\sim\mathsf{N}%
_{m}\left(  \boldsymbol{\mu},\boldsymbol{\Sigma}\right)  $ with
$\boldsymbol{\Sigma=\Sigma}_{\boldsymbol{\eta}}+\boldsymbol{\Sigma
}_{\mathbf{v}}$. On the other hand, for $\boldsymbol{\zeta}\sim\mathsf{N}%
_{m}\left(  \boldsymbol{\mu},\boldsymbol{\Sigma}\right)  $, using the spectral
decomposition of $\boldsymbol{\Sigma}$, we can construct uncorrelated
$\boldsymbol{\eta}\sim\mathsf{N}_{m}\left(  \mathbf{0},\boldsymbol{\Sigma
}_{\boldsymbol{\eta}}\right)  $ and $\mathbf{v}\sim\mathsf{N}_{m}\left(
\mathbf{0},\boldsymbol{\Sigma}_{\mathbf{v}}\right)  $ such that
\begin{equation}
\boldsymbol{\zeta}
=\boldsymbol{\mu}+\boldsymbol{\eta}+\mathbf{v}
\text{\quad and \quad}  \boldsymbol{\Sigma}={\Sigma}_{\boldsymbol{\eta}}+\boldsymbol{\Sigma}_{\mathbf{v}}. \label{eqModel}%
\end{equation}

We remark that the identifiability issue associated with the decomposition \eqref{eqModel}
has no effect on the validity of the SLLN for the sequence of (conditional) rejections or FDP
associated with the MTPs to be discussed next since the asymptotic analysis is determined by $\boldsymbol{\Sigma}$ and $\boldsymbol{\Sigma}_{\mathbf{v}}$. In particular, PFA in \cite{Fan:2012} and \cite{Fan:2017}
and the models in \cite{desai2012cross} and \cite{Friguet:2009} all employ such a decomposition. Further, the decomposition \eqref{eqModel} holds for decomposable multivariate probabilities including the Normal distribution; see \cite{Cuppens:1975} for more examples.

Write $\boldsymbol{\Sigma}=\left(  \tilde{\sigma}%
_{ij}\right)  $, let $\Phi$ be the CDF of $\mathsf{N}_{1}\left(  0,1\right)
$, and denote the standard deviation and variance of $v_{i}$ by $\sigma_{i,m}$ and $\sigma_{i,m}^{2}$.
Throughout the rest of the article, $\boldsymbol{\Sigma}$ is not necessarily a correlation matrix or invertible but we will exclude the trivial situation where the variance $\tilde{\sigma}_{ii}$ of $\zeta_i$ is $0$ for some $1 \le i \le m$.

\subsection{Marginal multiple testing}
\label{secMarginalTesting}

Recall $\boldsymbol{\zeta}$ has mean $\boldsymbol{\mu}=\left(  \mu
_{1},\ldots,\mu_{m}\right)  ^{\top}$. Consider multiple testing the $i$th null
hypothesis $H_{i0}:\mu_{i}=0$ versus $H_{i1}:\mu_{i}\neq0$ for all $1\leq
i\leq m$. Let $Q_{0,m}$ be the set of indices of the true null hypotheses
whose cardinality $\left\vert Q_{0,m}\right\vert $ is $m_{0}$, $Q_{1,m}$ that
for the false null hypotheses, and $\pi_{0,m}=m^{-1}m_{0}$ the
proportion of true null hypotheses. Given an observation
$\boldsymbol{\zeta}=\left(  \zeta_{1},\ldots,\zeta_{m}\right)  ^{\top}$,
define $p_{i}=1-F_{i}\left(  \zeta_{i}  \right)  $ as the
one-sided p-value and $p_{i}=2F_{i}\left(  -\left\vert \zeta_{i}\right\vert
\right)  $ as the two-sided p-value for $\zeta_{i}$, where $F_{i}$ is the
CDF of $\zeta_{i}$ when $\mu_{i}=0$.
We refer to these p-values as ``marginal p-values''.

Consider the marginal MTP with a rejection threshold $t\in\left(  0,1\right)  $ that
rejects $H_{i0}$ if and only if (iff) $p_{i}\leq t$. Then it induces
$R_{m}\left(  t\right)  =\sum_{i=1}^{m}1_{\left\{  p_{i}\leq t\right\}  }$ as
the number of rejections and $V_{m}\left(  t\right)  =\sum_{i\in Q_{0,m}%
}1_{\left\{  p_{i}\leq t\right\}  }$ as the number of false discoveries, where
$1_{A}$ is the indicator of a set $A$. Further, the FDP and FDR of the MTP are
respectively
\begin{equation*}
\mathrm{FDP}_{m}\left(  t\right)  =\frac{V_{m}\left(  t\right)  }{R_{m}\left(
t\right)  \vee1}\text{ \ \ \ and \ \ \ \ }\mathrm{FDR}_{m}\left(  t\right)
=\mathbb{E}\left[  \mathrm{FDP}_{m}\left(  t\right)  \right]  , 
\end{equation*}
where $a\vee b=\max\left\{  a,b\right\}  $ and $\mathbb{E}$ denotes expectation with respect to the probability measure $\mbb{P}$.

\subsection{Conditional multiple testing}

\label{SecCondMTP}

When $\boldsymbol{\Sigma}$ encodes strong dependence among the components
$\zeta_{i}$ of $\boldsymbol{\zeta}$, it is usually hard
to well estimate the FDP or FDR of the marginal MTP. However,
when $\boldsymbol{\eta}$ represents a dominant part of dependence among the
$\zeta_{i}$'s, the conditional p-values $\left\{p_i\vert \boldsymbol{\eta}\right\}_{i=1}^m$ have a much weaker dependence structure,
and a conditional MTP that assesses which $\mu_{i}$'s are $0$ based on $\left\{p_i\vert \boldsymbol{\eta}\right\}_{i=1}^m$ may have a more stable FDP and provide valuable information on the FDP and FDR of the marginal MTP. Specifically, with a
rejection threshold $t\in\left(  0,1\right)  $, the conditional MTP rejects
$H_{i0}:\mu_{i}=0$ iff $\left(  p_{i}|\eta_{i}\right)  \leq t$, i.e., it
rejects $H_{i0}$ iff the p-value $p_{i}$ conditional on $\eta_{i}$ is no larger than $t$.
This strategy has been taken by \cite{Friguet:2009} and \cite{Fan:2012} to study the performance of the marginal MTP
related to the Normal means problem, and will be referred to as the ``conditional MTP''.

In contrast to the previous conditional MTP, another conditional MTP can be defined as follows.
First, conditional on $\ms{\eta}$, compute the one-sided p-value
$\tilde{p}_i = 1-\tilde{F}_i \left(\zeta_i - \eta_i \right)$ or two-sided p-value
$\tilde{p}_i = 2 \tilde{F}_i \left(- \vert \zeta_i - \eta_i \vert \right)$, where $\tilde{F}_i$ is the CDF of
$\zeta_i - \eta_i$ when $\mu_i=0$. Then, with a rejection threshold $t\in\left[  0,1\right]  $, reject $H_{i0}:\mu_{i}=0$ iff $\tilde{p_{i}}  \leq t$. This approach has also been taken by \cite{Friguet:2009} under the name of ``factor-adjusted MTP''  and by \cite{Fan:2012} ``dependence-adjusted procedure'' as an adaptive MTP, and will be referred to as the ``adjusted conditional MTP''.
The adjusted conditional MTP is used more frequently than the conditional MTP. 

For the conditional MTP, let $X_{i}=1_{\left\{  p_{i}\leq t|\boldsymbol{\eta}\right\}
}$ be the indicator of whether the marginal p-value $p_{i}$ conditional on $\boldsymbol{\eta}$ is
no larger than $t$. Then $X_{i}=1_{\left\{  p_{i}\leq t|\eta_{i}\right\}  }$.
In other words, for the conditional MTP, $X_{i}$ is the indicator of whether $H_{i0}$ is rejected conditional on $\boldsymbol{\eta}$.
We call $\left\{  X_{i}\right\}  _{i=1}^{m}$ the
\textquotedblleft sequence of conditional rejections". The conditional MTP
induces the following quantities: the number of conditional rejections
$R_{m}\left(  t|\boldsymbol{\eta}\right)  =\sum_{i=1}^{m}X_{i}$, the number of
conditional false discoveries $V_{m}\left(  t|\boldsymbol{\eta}\right)
=\sum_{i\in Q_{0,m}}X_{i}$, the conditional FDP
\[
\mathrm{FDP}_{m}\left(  t|\boldsymbol{\eta}\right)  =\frac{V_{m}\left(
t|\boldsymbol{\eta}\right)  }{R_{m}\left(  t|\boldsymbol{\eta}\right)  \vee
1},
\]
and conditional FDR\ $\mathrm{FDR}_{m}\left(  t|\boldsymbol{\eta}\right)
=\mathbb{E}_{\mathbf{v}}\left[  \mathrm{FDP}_{m}\left(  t|\boldsymbol{\eta
}\right)  \right]  $. Here a random vector as a subscript of the expectation
$\mathbb{E}$ denotes the expectation with respect to the distribution of the
random vector.

\subsection{Principal correlation structure}
\label{SecDefPCS}

To quantify how much dependence among the $\zeta_{i}$'s
$\boldsymbol{\eta}$ should account for in the decomposition $\boldsymbol{\zeta
}=\boldsymbol{\mu}+\boldsymbol{\eta}+\mathbf{v}$, so that the conditional FDP
is well concentrated around its expectation, the conditional FDR, we introduce the concept of ``principal correlation
structure (PCS)". For a matrix $\mathbf{A}$ and $q>0$, let $\left\Vert
\mathbf{A}\right\Vert _{q}=\left(  \sum_{i,j}\left\vert \mathbf{A}\left(
i,j\right)  \right\vert ^{q}\right)  ^{1/q}$ be its $l_q$-norm. Let $\mathbf{R}_{\mathbf{v}}$ be the correlation matrix of $\mathbf{v}$. Define the
``covariance partition index'' for $\boldsymbol{\zeta
}$ in model (\ref{eqModel}) as $\varpi_{m}=m^{-2}\left\Vert \mathbf{R}_{\mathbf{v}}\right\Vert _{1}$. When $\varpi_{m}$ is small,
$\boldsymbol{\eta}$ captures the major part of the covariance dependence for
$\boldsymbol{\zeta}$ and $\mathbf{v}$ has less dependent components in terms of correlation. In
other words, when $\varpi_{m}$ is small, the conditional FDP may concentrate
around its expectation. When $\varpi_{m}$ is suitably small such that%
\begin{equation}
\varpi_{m}=m^{-2}\left\Vert \mathbf{R}_{\mathbf{v}}\right\Vert
_{1}=O\left(  m^{-\delta}\right)  \label{eqPCSdef}%
\end{equation}
for some $\delta>0$, where $O\left(  \cdot\right)  $ denotes Landau's big O
notation, we say that $\boldsymbol{\zeta}$ in (\ref{eqModel}) has a PCS.
Further, we call $\boldsymbol{\eta}$ the
\textquotedblleft principal vector" and $\mathbf{v}$ the \textquotedblleft
minor vector". Note that PCS restricts the order of the $l_1$-norm of the correlation matrix instead of the covariance
matrix of $\mb{v}$

We remark that $\varpi_{m}$ measures the relative magnitude of $\left\Vert \mathbf{R}_{\mathbf{v}}\right\Vert _{1}$ to $m^{2}$.
So, in the definition of PCS in \eqref{eqPCSdef}, $\varpi_{m}$
asymptotically proportional to $m^{-\delta}$ and decreasing should be interpreted as \textquotedblleft the speed
of increase of $\left\Vert \mathbf{R}_{\mathbf{v}}\right\Vert _{1}$ is slower in
order than $m^{2}$\textquotedblright. In other words, the absolute correlations among
components of the minor vector $\mathbf{v}$ accumulate at a rate slower than
the square of its dimension.

\section{SLLN for Normal means problem under dependence}

\label{secNormalMeans}

In this section, we derive our key results on Normal means problem under dependence. We will start with Hermite polynomials in \autoref{sec:Mehler}, and give in
\autoref{secVarianceOfCondRej} the exact formula for the variance of the average
number of conditional rejections $m^{-1}R_{m}\left(  t|\boldsymbol{\eta
}\right)  $ and an upper bound for this variance. Then we will show the SLLN for
the number of rejections and other random processes for the marginal MTP in \autoref{secSLLNmarginal} and for the conditional MTP in \autoref{secSLLNcondrej}.

We will deal with the adjusted conditional MTP in \autoref{secAcmtp}
since the conditions and techniques to derive the SLLN associated with it are almost identical to those for the conditional MTP.

\subsection{Hermite polynomial and Mehler expansion}

\label{sec:Mehler}

We state some basic facts on Hermite polynomials and Mehler expansion, which
will be used in the proofs of our key results. Let $\phi\left(  x\right)
=\left(  2\pi\right)  ^{-1/2}\exp\left(  -x^{2}/2\right)  $, i.e., $\phi$ is
the standard Normal density, and $f_{\rho}$ be the density of standard
bivariate Normal random vector with correlation $\rho\in\left(  -1,1\right)
$, i.e.,%
\[
f_{\rho}\left(  x,y\right)  =\frac{1}{2\pi\sqrt{1-\rho^{2}}}\exp\left(
-\frac{x^{2}+y^{2}-2\rho xy}{2\left(  1-\rho^{2}\right)  }\right)  .
\]
Let
\begin{equation*}
H_{n}\left(  x\right)  =\left(  -1\right)  ^{n}\frac{1}{\phi\left(
x\right)  }\frac{d^{n}}{dx^{n}}\phi\left(  x\right)
\end{equation*}
be the $n$th Hermite
polynomial; see \cite{Feller:1971B} for such a definition of $H_{n}$. Then
Mehler's expansion in \cite{Mehler:1866} implies%
\begin{equation}
f_{\rho}\left(  x,y\right)  =\left(  1+\sum_{n=1}^{\infty}\frac{\rho^{n}}%
{n!}H_{n}\left(  x\right)  H_{n}\left(  y\right)  \right)  \phi\left(
x\right)  \phi\left(  y\right)  . \label{eq:Mehler}%
\end{equation}
By \cite{Watson:1933}, the series on the right hand side of (\ref{eq:Mehler})
as a trivariate function of $\left(  x,y,\rho\right)  $ is uniformly
convergent on each compact set of $\mathbb{R}\times\mathbb{R}\times\left(
-1,1\right)  $.

We now provide a very important bound on these polynomials:
\begin{lemma}\label{BndHermitePoly}
  For the Hermite polynomials $H_n(\cdot)$, there is some constant $K_{0}>0$ independent of $n$ and $y$ such that
\begin{equation}
\left\vert e^{-y^{2}/2}H_{n}\left(  y\right)  \right\vert \leq K_{0}\sqrt
{n!}n^{-1/12}e^{-y^{2}/4}\text{ \ for any\ }y\in\mathbb{R}.
\label{eqBoundHermite}%
\end{equation}
\end{lemma}

\subsection{Variance of the average number of conditional rejections}
\label{secVarianceOfCondRej}
Let $\mathbb{V}_{\cdot}$ with a subscript denote the variance with respect to
the distribution of the random vector in the subscript, and so do the
subscript in the covariance operator $\mathrm{cov}_{\cdot}$. Recall the one-sided p-value $p_{i}=1-\Phi\left( \tilde{\sigma}_{ii}^{-1/2} \zeta_{i}\right)$, two-sided p-value $p_{i}=2\Phi\left(  -\tilde{\sigma}_{ii}^{-1/2}\left\vert \zeta_{i}\right\vert \right)$, $X_{i}=1_{\left\{  p_{i}\leq t|\boldsymbol{\eta}\right\}  }$ and $R_{m}\left(
t|\boldsymbol{\eta}\right)  =\sum_{i=1}^{m}X_{i}$. To derive a formula for the
variance $\mathbb{V}_{\mathbf{v}}\left[  m^{-1}R_{m}\left(  t|\boldsymbol{\eta
}\right)  \right]  $, we
introduce some notations. For a one-sided p-value $p_{i}$, define $\tilde
{t}=\tilde{\sigma}_{ii}^{1/2}\Phi^{-1}\left(  1-t\right)  $, $r_{1,i}=\tilde{t}-\mu_{i}-\eta_{i}$ and
$r_{2,i}=-\infty$; for a two-sided p-value $p_{i}$, define $\tilde{t}%
=-\tilde{\sigma}_{ii}^{1/2}\Phi^{-1}\left(  2^{-1}t\right)  $, $r_{1,i}=\tilde{t}-\mu_{i}-\eta_{i}$ and
$r_{2,i}=-\tilde{t}-\mu_{i}-\eta_{i}$. Further, set $c_{l,i}=\sigma_{i,m}%
^{-1}r_{l,i}$ for $l=1,2$, let $\rho_{ij}$ be the correlation between $v_{i}$
and $v_{j}$ for $i\neq j$, and define%
\begin{equation}
\left\{
\begin{array}
[c]{c}%
E_{1,m}=\left\{  \left(  i,j\right)  :1\leq i,j\leq m,i\neq j,\left\vert
\rho_{ij}\right\vert <1\right\}  ,\\
E_{2,m}=\left\{  \left(  i,j\right)  :1\leq i,j\leq m,i\neq j,\left\vert
\rho_{ij}\right\vert =1\right\}  .
\end{array}
\right.  \label{eqSets}%
\end{equation}
Namely, $E_{2,m}$ records pairs $\left(  v_{i},v_{j}\right)  $ with $i\neq j$
such that $v_{i}$ and $v_{j}$ are linearly dependent almost surely (a.s.).

\begin{lemma}
\label{LmVarianceNormal}Set
\begin{equation}
I_{1}=m^{-2}\sum_{i=1}^{m}\mathbb{V}_{\mathbf{v}}\left[  X_{i}\right]
+m^{-2}\sum_{\left(  i,j\right)  \in E_{2,m}}\mathrm{cov}_{\mathbf{v}}\left(
X_{i},X_{j}\right)  \label{eqVarPa}%
\end{equation}
and $d_{n}\left(  c,c^{\prime}\right)  =H_{n}\left(  c\right)  \phi\left(
c\right)  -H_{n}\left(  c^{\prime}\right)  \phi\left(  c^{\prime}\right)  $
for $c,c^{\prime}\in\mathbb{R}$. Then%
\begin{equation}
\mathbb{V}_{\mathbf{v}}\left[  m^{-1}R_{m}\left(  t|\boldsymbol{\eta}\right)
\right]  =I_{1}+m^{-2}\sum_{\left(  i,j\right)  \in E_{1,m}}\sum_{n=1}%
^{\infty}\frac{\rho_{ij}^{n}}{n!}H_{n-1}\left(  c_{1,i}\right)  H_{n-1}\left(
c_{1,j}\right)  \phi\left(  c_{1,i}\right)  \phi\left(  c_{1,j}\right)
\label{eqVarPB}%
\end{equation}
for one-sided p-values, and%
\begin{equation}
\mathbb{V}_{\mathbf{v}}\left[  m^{-1}R_{m}\left(  t|\boldsymbol{\eta}\right)
\right]  =I_{1}+m^{-2}\sum_{\left(  i,j\right)  \in E_{1,m}}\sum_{n=1}%
^{\infty}\frac{\rho_{ij}^{n}}{n!}d_{n-1}\left(  c_{1,i},c_{2,i}\right)
d_{n-1}\left(  c_{1,j},c_{2,j}\right)  \label{eqVarPC}%
\end{equation}
for two-sided p-values.
\end{lemma}

\begin{proof}
Expand $\mathbb{V}_{\mathbf{v}}\left[  m^{-1}R_{m}\left(  t|\boldsymbol{\eta
}\right)  \right]  $ into summands involving integrals, use Mehler's expansion in \autoref{sec:Mehler} for
$\rho_{ij}$ with $\left(  i,j\right)  \in E_{1,m}$ for the integrands in the double integrals, and observe $H_{n-1}%
\left(  x\right)  \phi\left(  x\right)  =\int_{-\infty}^{x}H_{n}\left(
y\right)  \phi\left(  y\right)  dy$ for $x\in\mathbb{R}$, we get the results.
This completes the proof.
\end{proof}

\autoref{LmVarianceNormal} gives the exact value for the variance of
$m^{-1}R_{m}\left(  t|\boldsymbol{\eta}\right)  $. In case the SLLN fails, \autoref{LmVarianceNormal} can be
used, e.g., in combination with Markov inequality, to give a bound on the
deviation of $R_{m}\left(  t|\boldsymbol{\eta}\right)  $ from its mean
$\mathbb{E}_{\mathbf{v}}\left[  R_{m}\left(  t|\boldsymbol{\eta}\right)
\right]  $ for any $m\geq1$.

To obtain bounds on the variance $\mathbb{V}_{\mathbf{v}}\left[  m^{-1}%
R_{m}\left(  t|\boldsymbol{\eta}\right)  \right]  $, we introduce sets that
describe different behavior of the $\eta_{i}$'s or $v_{i}$'s. Define%
\begin{equation*}
E_{0}=\left\{  i\in\mathbb{N}:\sigma_{i,m}>0\text{ for any }m\text{ but
}\liminf_{m\rightarrow\infty}\sigma_{i,m}=0\right\}  
\end{equation*}
and set $E_{0,m}=E_{0}\cap\left\{  1,\ldots,m\right\}  $, i.e., $E_{0,m}$
contains $i$ such that the standard deviation $\sigma_{i,m}\ $of $v_{i}$ can
be arbitrarily small as $m\rightarrow\infty$. Further, define%
\begin{equation}
G_{m,\boldsymbol{\eta}}\left(  t,\varepsilon_{m}\right)  =\bigcup
\nolimits_{i\in E_{0,m}}\left\{  \omega\in\Omega:\min\left\{
\left\vert r_{1,i}\right\vert ,\left\vert r_{2,i}\right\vert \right\}
<\varepsilon_{m}\right\}  \label{eqExceptionB}%
\end{equation}
for some $\varepsilon_{m}>0$ (be determined later) such that $\lim
_{m\rightarrow\infty}\varepsilon_{m}=0$. Namely, $G_{m,\boldsymbol{\eta}%
}\left(  t,\varepsilon_{m}\right)  $ contains $\eta_{i}$ that is within
distance $\varepsilon_{m}$ from $\pm\tilde{t}-\mu_{i}$ and whose variance
 is arbitrarily close to that of $\zeta_{i}$ as
$m\rightarrow\infty$. Note that $G_{m,\boldsymbol{\eta}}\left(  t,\varepsilon
_{m}\right)  =\varnothing$ when $E_{0}=\varnothing$ and that the Cartesian
product $E_{0,m}\times E_{0,m}$ contains distinct $i$ and $j$ for which the
covariance
\[
\mathrm{cov}_{\mathbf{v}}\left(  X_{i},X_{j}\right)  =\mathbb{E}_{\mathbf{v}%
}\left[  X_{i}X_{j}\right]  -\mathbb{E}_{\mathbf{v}}\left[  X_{i}\right]
\mathbb{E}_{\mathbf{v}}\left[  X_{j}\right]
\]
may inflate the order of $\mathbb{V}_{\mathbf{v}}\left[  m^{-1}R_{m}\left(
t|\boldsymbol{\eta}\right)  \right]  $.

With the above preparations, we have:
\begin{proposition}
\label{ThmSLLN}There is a constant $C>0$ such that
\begin{equation}\label{IneqComp}
  \left\vert\mathrm{cov}_{\mathbf{v}}\left(  X_{i},X_{j}\right)\right\vert  \le C \left\vert \rho_{ij}\right\vert
\end{equation}
for any $1 \le i \le j \le m$, $m \ge 1$, $t \in \left(0,1\right)$, $\ms{\mu}$ and $\ms{\eta}$. Suppose for some $\delta>0$%
\begin{equation}
m^{-2}\left\Vert \boldsymbol{\Sigma}_{\mathbf{v}}\right\Vert _{1}=O\left(
m^{-\delta}\right)  \text{ \ and \ }\left\vert E_{2,m}\right\vert =O\left(
m^{2-\delta}\right)  . \label{eqKeyCondA}%
\end{equation}
Let
\begin{equation}\label{minStdMinor}
  \sigma_{0}=\lim_{m\rightarrow\infty}\min\left\{  \sigma_{i,m}:\sigma_{i,m}\neq0,1\leq i\leq m\right\}
\end{equation}
and $D_{m,\boldsymbol{\eta}}\left(  t,\varepsilon_{m}\right)  =\Omega\setminus
G_{m,\boldsymbol{\eta}}\left(  t,\varepsilon_{m}\right)  $. If $\sigma_{0}>0$,
then%
\begin{equation}
\mathbb{V}_{\mathbf{v}}\left[  m^{-1}R_{m}\left(  t|\boldsymbol{\eta
}  \right)  \right]  =O\left(  m^{-\min\left\{
\delta,1\right\}  }\right)  ; \label{eqKeyVarBound}%
\end{equation}
if $\sigma_{0}=0$, then for each $\ms{\eta}\rbk{\omega}$ with $\omega \in D_{m,\boldsymbol{\eta}}\left(t,\varepsilon_{m}\right)$,%
\begin{equation}
\mathbb{V}_{\mathbf{v}}\left[  m^{-1}R_{m}\left(  t|\ms{\eta}
  \right)  \right]  =O\left(  \varepsilon
_{m}^{-2}m^{-\min\left\{  \delta,1\right\}  }\right)  . \label{eqKeyVarBoundA}%
\end{equation}

\end{proposition}

The inequality (\ref{IneqComp}) is the universal comparison result that connects the correlation $\rho_{ij}$ between $v_i$ and $v_j$ with the covariance between
the indicator functions of the p-values of their associated tests. When $\boldsymbol{\eta} = \mathbf{0}$ a.s., (\ref{IneqComp}) connects the correlation between the two components of a bivariate Normal random vector to the covariance of the indicator functions of the p-values of testing each of the two marginal means of the vector. Further, with a bit more effect on studying the asymptotics of Hermite polynomials related to \autoref{BndHermitePoly}, the constant $C$ in (\ref{IneqComp}) can be identified.
\autoref{ThmSLLN} implies that, when there are not excessively many linearly
dependent pairs $\left(  v_{i},v_{j}\right)  $, $i\neq j$ and
$\boldsymbol{\zeta}$ has a PCS, $\mathbb{V}_{\mathbf{v}}\left[  m^{-1}%
R_{m}\left(  t|\boldsymbol{\eta}\right)  \right]  $ is of order $m^{-\delta}$
when the limit $\sigma_{0}$ of the minimum of the nonzero standard deviations
$\sigma_{i,m}$ for the $v_{i}$'s is positive, whereas $\mathbb{V}_{\mathbf{v}%
}\left[  m^{-1}R_{m}\left(  t|\ms{\eta}\right)  \right]  $ is
of order$\ m^{-\delta}\varepsilon_{m}^{-2}$ for
$\ms{\eta}\rbk{\omega}$ with $\omega \in D_{m,\boldsymbol{\eta}}\left(t,\varepsilon_{m}\right)$ if $\sigma_{0}=0$.

\subsection{SLLN for the marginal multiple testing procedure}
\label{secSLLNmarginal}

We deal with the marginal MTP where $\left\{\zeta_i, 1 \le i \le m, m \ge 1\right\}$ does not form a triangular array. In terms of the additive decomposition \eqref{eqModel}, this setting is equivalent to $\ms{\eta} = \mathbf{0}$ a.s. for all $m \ge 1$ and the sequence $\left\{v_i, 1 \le i \le m, m \ge 1\right\}$ not being a triangular array, which implies $\ms{\Sigma}_{\ms{\eta}}=\mb{0}$, $\ms{\Sigma}=\ms{\Sigma}_{\mb{v}}$, $\ms{\zeta}=\ms{\mu}+\mb{v}$ a.s. for all $m \ge 1$. Further, the conditional MTP is the marginal MTP. Our first result characterizes the type of dependence via PCS for which the SLLN holds for the marginal MTP:

\begin{theorem}\label{prop:SLLNMarginal}
  Assume $\boldsymbol{\zeta}\sim\mathsf{N}_{m}\left(  \boldsymbol{\mu},\boldsymbol{\Sigma}\right)  $ with correlation matrix $\mb{R}$.
  If
  \begin{equation}
m^{-2}\left\Vert \mathbf{R}\right\Vert _{1}=O\left(  m^{-\delta
}\right)  \text{ \ for \ some }\delta>0, \label{eqCondStrongest}%
\end{equation}
then
  $m^{-1}\left\vert R_{m}\left(  t  \right)  -\mathbb{E}\left[  R_{m}\left(t  \right)  \right]\right\vert \to 0$ a.s. and
  $m^{-1}\left\vert V_{m}\left(  t  \right)  -\mathbb{E}\left[  V_{m}\left(t  \right)  \right]\right\vert \to 0$ a.s. as $m \to \infty$.
If further $\liminf_{m\rightarrow\infty}\pi_{0,m}  >0$, then
$\left\vert \mathrm{FDP}_{m}\left(  t\right)  -\mathbb{E}\left[
\mathrm{FDP}_{m}\left(  t  \right)  \right]  \right\vert \to 0$ a.s. as $m \to \infty$.
\end{theorem}

In short, when the correlation matrix $\mb{R}$ of $\ms{\zeta}$ satisfies $m^{-2}\left\Vert \mb{R}\right\Vert _{1}=O\left(m^{-\delta}\right)$ for some $\delta >0$, there is no need for the conditional MTP or the adjusted conditional MTP, and the marginal MTP has well concentrated sequences of rejections and FDP in terms of the SLLN. To apply \autoref{prop:SLLNMarginal}, we only need to verify the order of $m^{-2}\left\Vert \mathbf{R}\right\Vert _{1}$, and estimate the variances $\tilde{\sigma}_{ii}$ rather than $\mathbf{R}$ accurately.
In contrast, if we only restrict the absolute covariances between the components of $\ms{\zeta}$, then we may only have a partial SLLN that excludes certain pairs of $\rbk{\ms{\mu},\ms{\sigma}_{\mb{v}}}$, where $\ms{\sigma}_{\mb{v}}=\left(\sigma_{1,m},\ldots,\sigma_{m,m}\right)$:
\begin{theorem}
\label{ThmSLLNAMarg}Suppose $\ms{\eta} = \mb{0}$ a.s. for all $m \ge 1$ and that the sequence $\left\{\zeta_i, 1 \le i \le m, m \ge 1\right\}$ is not a triangular array.
Assume \eqref{eqKeyCondA}, i.e., $m^{-2}\left\Vert
\boldsymbol{\Sigma}_{\mathbf{v}}\right\Vert _{1}=O\left(  m^{-\delta}\right)
$ and\ $\left\vert E_{2,m}\right\vert =O\left(  m^{2-\delta}\right)  $ for some $\delta>0$.
 Set
$\varepsilon_{m}=m^{-\delta_{1}}$ for any $\delta_{1}\in\left(  0,\min\left\{
2^{-1}\delta,2^{-1}\right\}  \right)  $ and
\begin{equation}\label{eqSetImpossibleMarg}
G_{t}=\bigcup\nolimits_{m\geq1}\left\{ \mu_i,\tilde{\sigma}_{ii}:i\in
E_{0,m},\min\left\{  \left\vert r_{1,i}\right\vert ,\left\vert r_{2,i}%
\right\vert \right\}  <\varepsilon_{m}\right\}.
\end{equation}
Then  $m^{-1}\left\vert R_{m}\left(  t  \right)  -\mathbb{E}\left[  R_{m}\left(t  \right)  \right]\right\vert \to 0$ a.s. and
  $m^{-1}\left\vert V_{m}\left(  t  \right)  -\mathbb{E}\left[  V_{m}\left(t  \right)  \right]\right\vert \to 0$ a.s. as $m \to \infty$ when either (i) $\sigma_0 >0$ or (ii) $\sigma_{0} =0 $ and $\rbk{\ms{\mu},\ms{\sigma}_{\mb{v}}}\notin G\left(  t\right)$. If in addition $\liminf_{m \to \infty}\pi_{0,m}>0$, then $\left\vert \mathrm{FDP}_{m}\left(  t\right)  -\mathbb{E}\left[
\mathrm{FDP}_{m}\left(  t  \right)  \right]  \right\vert \to 0$ a.s. as $m \to \infty$ when either (i) $\sigma_0 >0$ or (ii) $\sigma_{0} =0 $ and $\rbk{\ms{\mu},\ms{\sigma}_{\mb{v}}}\notin G\left(  t\right)$.
\end{theorem}

Note that $\sigma_0$ is defined by (\ref{minStdMinor}) and that $G_t$ in \eqref{eqSetImpossibleMarg} restricts the joint behavior of
the mean vector $\ms{\mu}$ and the vector of variances $\ms{\sigma}_{\ms{\zeta}}=\rbk{\tilde{\sigma}_{11},\cdots,\tilde{\sigma}_{mm}}^{\top}$
of the Normal random vector $\ms{\zeta}$.
Comparing \autoref{ThmSLLNAMarg} and \autoref{prop:SLLNMarginal}, we see that, for the marginal MTP
\begin{equation}\label{pfaKey}
  m^{-2}\left\Vert \boldsymbol{\Sigma}\right\Vert _{1}=O\left(m^{-\delta}\right)
\end{equation}
may be insufficient to ensure the SLLN for the sequence of rejections, whereas (\ref{eqCondStrongest}), i.e., $m^{-2}\left\Vert \mb{R}\right\Vert _{1}=O\left(m^{-\delta}\right)$ is. Condition (\ref{eqCondStrongest}) excludes
cases for which the covariance matrix of $\ms{\zeta}$ has a small magnitude
but the correlations among components of $\ms{\zeta}$ are still strong enough
to invalidate the SLLN for the marginal MTP. On the the hand, it is easy to see that $\sigma_0 >0$, $\lim_{m \to \infty}\sup_{1\le i \le m} \tilde{\sigma}_{ii} < \infty$ and (\ref{pfaKey}) together implies (\ref{eqCondStrongest}). In other words, the covariance matrix of $\ms{\zeta}$ usually cannot be singular asymptotically if (\ref{pfaKey}) were to induce the SLLN for the marginal MTP.

\subsection{SLLN for the conditional multiple testing procedure}
\label{secSLLNcondrej}
Recall the additive decomposition \eqref{eqModel}, i.e., $\boldsymbol{\zeta}=\boldsymbol{\mu}+\boldsymbol{\eta}+\mathbf{v}$.
When the major vector $\ms{\eta}$ is not $\mathbf{0}$ a.s. for all $m \ge 1$ and the components $v_i, 1 \le 1 \le m$ of $\mathbf{v}$ form a triangular array as $m$ changes, techniques for the SLLN of triangular arrays are needed. In the rest of the paper, whenever needed, a subscript $m$ will be added to a quantity to indicate its dependence on $m$; e.g., $\ms{\mu}$, $\ms{\eta}$ and $\mb{v}$ will also be written respectively as $\ms{\mu}_m$, $\ms{\eta}_m$ and $\mb{v}_m$.

For an $m^{\prime}$-dimensional vector $\mathbf{a}=\left(  a_{1},\ldots
,a_{m^{\prime}}\,\right)  $ and a natural number $m\leq m^{\prime}$, let
$\mathbf{a}^{\left(  m\right)  }=$ $\left(  a_{1},\ldots,a_{m}\,\right)  $.
Let $\ms{\sigma}_{\mb{v},m}=\left(\sigma_{1,m},\ldots,\sigma_{m,m}\right)$ be the vector of
standard deviations of $\mb{v}_m$ and recall $\sigma_0$ defined by (\ref{minStdMinor}). We have the following result:

\begin{theorem}
\label{ThmSLLNStrongest}
Assume $\sigma_0>0$ and
   \begin{equation}
m^{-2}\left\Vert \mathbf{R}_{\mathbf{v}}\right\Vert _{1} = O\left(   m^{-\delta
} \right) \text{ \ for \ some }\delta >0. \label{eqCondStrongestA}%
\end{equation}
If for any natural numbers
$m,m^{\prime}$ such that $m\leq m^{\prime}$ and $m\rightarrow\infty$%
\begin{equation}
\max\left\{  \left\Vert \boldsymbol{\eta}_{m}-\boldsymbol{\eta}_{m^{\prime}%
}^{\left(  m\right)  }\right\Vert _{2},\left\Vert \boldsymbol{\mu}%
_{m}-\boldsymbol{\mu}_{m^{\prime}}^{\left(  m\right)  }\right\Vert
_{2},\left\Vert \boldsymbol{\sigma}_{\mathbf{v},m}-\boldsymbol{\sigma
}_{\mathbf{v},m^{\prime}}^{\left(  m\right)  }\right\Vert _{2}\right\}
\rightarrow0\text{ a.s.,}\label{CondHomo}%
\end{equation}
then
$m^{-1}\left\vert R_{m}\left(  t|\ms{\eta}  \right)  -\mathbb{E}_{\mathbf{v}}\left[  R_{m}\left(
t|\ms{\eta}  \right)  \right]
\right\vert \to 0$ a.s and $m^{-1}\left\vert V_{m}\left(  t|\ms{\eta}  \right)  -\mathbb{E}_{\mathbf{v}}\left[  V_{m}\left(
t|\ms{\eta}  \right)  \right] \right\vert \to 0$ a.s.
If in addition $\liminf_{m\rightarrow\infty}\pi_{0,m}  >0$,  then
$m^{-1}\left\vert \mathrm{FDP}_{m}\left(  t|\ms{\eta}\right)
 -\mathbb{E}_{\mathbf{v}}\left[  \mathrm{FDP}_{m}\left(  t|\ms{\eta}\right)  \right]
\right\vert \to 0$ a.s.
\end{theorem}

\autoref{ThmSLLNStrongest} implies that $\boldsymbol{\zeta}$ having a PCS, i.e.,
$m^{-2}\left\Vert \mathbf{R}_{\mathbf{v}}\right\Vert _{1}=O\left(m^{-\delta}\right)$
for some $\delta>0$, is usually insufficient to ensure the SLLN for the conditional MTP.
In \autoref{SecExample}, we will show by an example that $\boldsymbol{\zeta}$ having a PCS is almost necessary for such a SLLN to
hold. The condition \eqref{CondHomo} is referred to as the ``homogeneity condition'' on the triangular arrays of components of  $\ms{\mu}_m$, $\ms{\eta}_m$ and $\mb{v}_m$.
If the decomposition $\boldsymbol{\zeta}=\boldsymbol{\mu}+\boldsymbol{\eta}+\mathbf{v}$ is obtained from the spectral decomposition of the covariance matrix $\ms{\Sigma}$ of $\ms{\zeta}$ and components of $\ms{\mu}$ does not form a triangular array, then the homogeneity condition requires that the spectral decomposition of $\ms{\Sigma}$ affects little $\mathbf{v}$ or $\ms{\eta}$ as $\ms{\Sigma}$ changes.

Condition \eqref{CondHomo} is not hard to check, seems to be restrictive and can perhaps be weakened without invalidating the SLLN. However, we will not pursue it here and point out that certain types of homogeneity for the sequences $\ms{\mu}_m$, $\ms{\eta}_m$ and $\mb{v}_m$ as $m$ changes are needed in order to obtain the SLLN.
On the other hand, condition (\ref{eqCondStrongestA}) induces the SLLN on a subsequence of the sequence of rejections, and the conditions $\sigma_0>0$ and \eqref{CondHomo} together with the continuity of Normal CDFs induce a controlled maximal inequality in the interpolation step of justifying the SLLN for triangular arrays.


\section{An example for which the SLLN fails without PCS}

\label{SecExample}

We provide an example for which the SLLN fails without PCS. Recall $\boldsymbol{\zeta}\sim\mathsf{N}_{m}\left(  \boldsymbol{\mu
},\boldsymbol{\Sigma}\right)  $.
Consider the representation $\boldsymbol{\zeta}=\boldsymbol{\mu}%
+\boldsymbol{\eta}+\mathbf{v}$ where $\boldsymbol{\eta}\sim\mathsf{N}%
_{m}\left(  \mathbf{0},\boldsymbol{\Sigma}_{\boldsymbol{\eta}}\right)  $ and
$\mathbf{v}\sim\mathsf{N}_{m}\left(  \mathbf{0},\boldsymbol{\Sigma
}_{\mathbf{v}}\right)  $ are uncorrelated and $\boldsymbol{\Sigma
}=\boldsymbol{\Sigma}_{\boldsymbol{\eta}}+\boldsymbol{\Sigma}_{\mathbf{v}%
}=\left(  \tilde{\sigma}_{ij}\right)  $.
Recall $X_{i}=1_{\left\{  p_{i}\leq
t|\boldsymbol{\eta}\right\}  }$ and $R_{m}\left(  t|\boldsymbol{\eta}\right)
=\sum_{i=1}^{m}X_{i}$. So, $\left\{  X_{i}:1\leq i\leq m\right\}  $ is a
sequence of dependent Bernoulli random variables, and $m^{-1}R_{m}\left(
t|\boldsymbol{\eta}\right)  $ is the average location of the \textquotedblleft
random walk" induced by $\left\{  X_{i}\right\}  _{i=1}^{m}$. We will
write $X_{i}\ $as $X_{i}\left(  t,\mb{v}\rbk{\omega}|\tilde{\ms{\eta}}
  \right)  $ when $\boldsymbol{\eta}\rbk{\omega}=\tilde{\ms{\eta}}$.
The following example illustrates that PCS is almost necessary for the
SLLN to hold.

\begin{proposition}
\label{cor:Failure}
For $m\geq3$ there exist a sequence of $\boldsymbol{\zeta}_{m}\sim
\mathsf{N}_{m}\left(  \mathbf{0},\boldsymbol{\Sigma}_{m}\right)  $ such that
$\boldsymbol{\zeta}_{m}=\boldsymbol{\eta}+\mathbf{v}$ for two uncorrelated
Normal random vectors $\boldsymbol{\eta}$ and $\mathbf{v}$. However, for this
sequence there exits a set $H_{t}\in\mathcal{F}$ with $\mathbb{P}\left(
H_{t}\right)  >0$ such that the SLLN fails for $\left\{  X_{i}\left(
t,\mb{v}\rbk{\omega}|\tilde{\ms{\eta}}  \right)  :1\leq
i\leq m, m \ge 1\right\}  $ for each $\tilde{\ms{\eta}}\in Q_t = \left\{\ms{\eta}\rbk{\omega}: \omega \in H_{t}\right\}$.
\end{proposition}

\begin{proof}
First, we construct the covariance matrices $\boldsymbol{\Sigma}%
_{\boldsymbol{\eta}}$ and $\boldsymbol{\Sigma}_{\mathbf{v}}$. Let
$\boldsymbol{\tilde{\gamma}}_{1}=\left(  \frac{-\sqrt{2}}{2},\frac{\sqrt{2}%
}{2},0,\ldots,0\right)  ^{\top}$, $\boldsymbol{\tilde{\gamma}}_{2}=\left(
\frac{\sqrt{2}}{2},\frac{-\sqrt{2}}{2},0,\ldots,0\right)  ^{\top}$ and
$\boldsymbol{\tilde{\gamma}}_{3}=\mathbf{1}_{m}$, where $\mathbf{1}_{m}$ is a
column of vector of $m$ $1$'s. Then $\boldsymbol{\tilde{\gamma}}_{i}%
^{\top}\boldsymbol{\tilde{\gamma}}_{j}=0$ when $i\neq j$. Let $\boldsymbol{\Sigma
}_{\boldsymbol{\eta}}=\boldsymbol{\tilde{\gamma}}_{1}\boldsymbol{\tilde
{\gamma}}_{1}^{\top}$, $\mathbf{\tilde{T}}=\left(  \boldsymbol{\tilde{\gamma}%
}_{3},\boldsymbol{\tilde{\gamma}}_{2}\right)  $ and $\boldsymbol{\Sigma
}_{\mathbf{v}}=\mathbf{\tilde{T}\tilde{T}}^{\top}$.

Secondly, we construct the sequence of Normal random vectors $\left\{
\boldsymbol{\zeta}_{m}\right\}  _{m}$, each with decomposition
$\boldsymbol{\zeta}_{m}=\boldsymbol{\eta}+\mathbf{v}$ for two uncorrelated
Normal random vectors $\boldsymbol{\eta}$ and $\mathbf{v}$. Let $w_{1}%
\sim\mathsf{N}_{1}\left(  0,1\right)  $ and $\mathbf{\tilde{w}}_{2}=\left(
w_{2},w_{3}\right)  ^{\top}\sim\mathsf{N}_{2}\left(  \mathbf{0},\mathbf{I}%
_{2}\right)  $ such that $w_{1}$ and $\mathbf{\tilde{w}}_{2}$ are independent.
Set $\boldsymbol{\eta}=\boldsymbol{\tilde{\gamma}}_{1}w_{1}$ and
$\mathbf{v}=\mathbf{\tilde{T}\tilde{w}}_{2}$. Then $\boldsymbol{\eta}%
\sim\mathsf{N}_{m}\left(  \mathbf{0},\boldsymbol{\Sigma}_{\boldsymbol{\eta}%
}\right)  $ and $\mathbf{v}\sim\mathsf{N}_{m}\left(  \mathbf{0}%
,\boldsymbol{\Sigma}_{\mathbf{v}}\right)  $, and $\boldsymbol{\eta}$ is
uncorrelated with $\mathbf{v}$. Note that $\boldsymbol{\Sigma}%
_{\boldsymbol{\eta}}$ and $\boldsymbol{\Sigma}_{\mathbf{v}}$ are singular. Set
$\boldsymbol{\zeta}_{m}=\boldsymbol{\eta}+\mathbf{v}$. Then $\boldsymbol{\zeta
}_{m}\sim\mathsf{N}_{m}\left(  \mathbf{0},\boldsymbol{\Sigma}_{m}\right)  $
and $\boldsymbol{\Sigma}_{m}=\boldsymbol{\Sigma}_{\boldsymbol{\eta}%
}+\boldsymbol{\Sigma}_{\mathbf{v}}$. Note that $\boldsymbol{\Sigma}_{m}$ is singular since
$\mathrm{rank}\left(  \Sigma_{m}\right)  \leq3$. Let $\boldsymbol{\eta
}=\left(  \eta_{1},\ldots,\eta_{m}\right)  ^{\top}$ and $\mathbf{v}=\left(
v_{1},\ldots,v_{m}\right)  ^{\top}$. Then,
\begin{equation}
\eta_{1}=-\frac{\sqrt{2}}{2}w_{1},\text{ }\eta_{2}=\frac{\sqrt{2}}{2}%
w_{1}\text{ \ and \ }\eta_{i}=0\text{ \ for \ }3\leq i\leq m, \label{eqEgA}%
\end{equation}
and
\begin{equation}
v_{1}=w_{2}+\frac{\sqrt{2}}{2}w_{3},\text{ }v_{2}=w_{2}-\frac{\sqrt{2}}%
{2}w_{3}\text{ \ and \ }v_{i}=w_{2}\text{ \ for \ }3\leq i\leq m.
\label{eqEgB}%
\end{equation}

Finally, we show that the SLLN fails for $\left\{  X_{i}:1\leq i\leq m,
m \ge 1\right\}  $. Recall $\tilde{t}=-\Phi^{-1}\left(  2^{-1}t\right)  $,
$r_{1,i}=\tilde{t}-\eta_{i}$, $r_{2,i}=-\tilde{t}-\eta_{i}$ for two-sided
p-values or $r_{2,i}=-\infty$ for one-sided p-values, and $c_{l,i}%
=\sigma_{i,m}^{-1}r_{l,i}$ for $l=1,2$. Define%
\[
A_{i}=\left\{  r_{2,i}\leq v_{i}\leq r_{1,i}\right\}
\]
for $1\leq i\leq m$. Then $A_{1}=\left\{  r_{2,1}\leq v_{1}\leq r_{1,1}%
\right\}  $ and $A_{2}=\left\{  r_{2,2}\leq v_{2}\leq r_{1,2}\right\}  $.
Further, for\ $3\leq i\leq m$, $A_{i}=\left\{  -\tilde{t}\leq w_{2}\leq
\tilde{t}\right\}  $\ for two-sided p-values and $A_{i}=\left\{  -\infty\leq
w_{2}\leq\tilde{t}\right\}  $ for one-sided p-values.

Let $Y_{i}=1_{\left\{  p_{i}\geq t|\boldsymbol{\eta}\right\}  }$, $\theta
_{i}=\mathbb{E}_{\mathbf{v}}\left[  Y_{i}\right]  $, $\bar{Y}_{m}=m^{-1}%
\sum_{i=1}^{m}Y_{i}$ and $\bar{\theta}_{m}=m^{-1}\sum_{i=1}^{m}\theta_{i}$.
Since $\left(  p_{i}|\boldsymbol{\eta}\right)  \geq t$ iff $\left\vert
\zeta_{i}\right\vert \leq\tilde{t}$ iff $v_{i}\in A_{i}$, we have%
\[
\theta_{i}=\int_{A_{i}}\frac{1}{\sqrt{2\pi}}\exp\left(  -\frac{1}{2}%
x^{2}\right)  dx=\int_{c_{2,i}}^{c_{1,i}}\frac{1}{\sqrt{2\pi}}\exp\left(
-\frac{1}{2}x^{2}\right)  dx
\]
and%
\begin{equation}
\mathbb{P}\left(  \bar{Y}_{m}-\bar{\theta}_{m}=1-\bar{\theta}_{m}\right)
=\mathbb{P}\left(  v_{1}\in A_{1},v_{2}\in A_{2},v_{1}\in A_{3}\right)
\label{eqProb}%
\end{equation}
conditional on $\boldsymbol{\eta}$. Clearly, there exits a set $H_{t}%
\in\mathcal{F}$ independent of $m$ such that: (i) $\mathbb{P}\left(
H_{t}\right)  >0$, (ii) $\limsup_{m\rightarrow\infty}\max_{1\leq i\leq
m}\theta_{i}<1$ conditional on each $\tilde{\ms{\eta}}\in Q_t$,
where $Q_t = \left\{\ms{\eta}\rbk{\omega}: \omega \in H_{t}\right\}$, and (iii) the right hand side of
(\ref{eqProb}) is positive conditional on $\tilde{\ms{\eta}}\in Q_t$. Thus, conditional on
$\tilde{\ms{\eta}}\in Q_t$,
\begin{equation}
\mathbb{P}\left(  \limsup_{m\rightarrow\infty}\left\vert 1-\bar{\theta}%
_{m}\right\vert >0\right)  >0. \label{eqFailure}%
\end{equation}
Since%
\[
-\left(  \bar{Y}_{m}-\bar{\theta}_{m}\right)  =m^{-1}R_{m}\left(
t,\mb{v}\rbk{\omega}|\tilde{\ms{\eta}}  \right)
-\mathbb{E}_{\mathbf{v}}\left[  m^{-1}R_{m}\left(  t,\mb{v}\rbk{\omega}|\tilde{\ms{\eta}}  \right)  \right]  ,
\]
(\ref{eqFailure}) implies that the SLLN does not hold for $\left\{
X_{i}:1\leq i\leq m, m \ge 1\right\}  $. This completes the proof.
\end{proof}

In the example provided by \autoref{cor:Failure}, the failure of the SLLN for
$\left\{  X_{i}:1\leq i\leq m, m \ge 1\right\}  $ is mainly due to $m^{-2}\left\Vert
\boldsymbol{\Sigma}_{\mathbf{v}}\right\Vert _{1}=O\left(  1\right)  $ and that
there are $O\left(  m^{2}\right)  $ linearly dependent pairs $\left(
v_{i},v_{j}\right)  $, $i\neq j$. In this case, $\boldsymbol{\zeta}_{m}$ does
not have a PCS, and $\left\{  m^{-1}R_{m}\left(  t|\boldsymbol{\eta}\right)
:m\geq1\right\}  $ is dominated by a random walk induced by components of
$\mathbf{v}$ and $\boldsymbol{\eta}$ given by (\ref{eqEgA}) and (\ref{eqEgB}).

\section{Simulation study}
\label{secSim}

We present a simulation study to verify our theory and compare PCS with PFA.
The validity of the SLLN for the sequence of (conditional) rejections is assessed by checking the sample variance of the average number of (conditional) rejections obtained from a large number of i.i.d. experiments. If the sequence of sample variances indexed by the dimension $m$ of the Normal random vector displays a strong trend of converging to $0$ as $m$ becomes larger and larger, we accept the validity of the SLLN; otherwise, we reject it. For example, if the sample variance of $m^{-1}R_{m}\left(  t|\boldsymbol{\eta}\right), m \ge 1$ shows a strong trend of converging to $0$ as $m$ increases, we accept that the SLLN holds for $R_{m}\left(  t|\boldsymbol{\eta}\right), m \ge 1$.

The simulation design in given in \autoref{simDesign}. We set the rejection threshold $t$ so that $m^{-1}R_{m}\left(  t|\boldsymbol{\eta}\right)$ is not constant for all independent repetitions of each simulation setting with a fixed combination of dependence type, sparsity regime and value of $m$, and that it is not zero for a repetition of the sequence of simulation settings with the same dependence type and sparsity regime but increasing $m$.
This helps prevent any spurious convergence to $0$ of a sequence of sample variances indexed by $m$ due to its corresponding sequence of average numbers of (conditional) rejections being constant (or $0$ in particular).
When we report the simulation results, we simply refer to ``sample variance'' as ``variance''.

\subsection{Simulation design}\label{simDesign}

Recall that $\boldsymbol{\zeta}\sim\mathsf{N}_{m}\left(  \boldsymbol{\mu
},\boldsymbol{\Sigma}\right)  $ with $\boldsymbol{\mu}=\left(  \mu_{1},\ldots,\mu_{m}\right)  ^{\top}$.
We consider $7$ values for $m$ as $500$, $1000$, $2000$, $4000$, $6000$, $8000$ or $10000$.
In order to compare our method with PFA in
\cite{Fan:2012}, we set $\boldsymbol{\Sigma}=\left(\tilde{\sigma}_{ij}\right)$ as a correlation matrix and consider
$6$ types of dependence structure encoded by $\boldsymbol{\Sigma}$.
Recall $m_{0}$ as the number of zero $\mu_{i}$'s and set $\pi_{1,m}=1-m_{0}m^{-1}$ as the
proportion of nonzero $\mu_{i}$'s.
We consider $3$ sparsity regimes, i.e., $\pi_{1,m}=0.05,m^{-0.4}$ or
$m^{-0.7}$, corresponding to the dense, moderately sparse, and very sparse regime (termed so as in \cite{Jin:2008}).
The nonzero $\mu_{i}$'s are generated independently such that their absolute values $\vert \mu_i \vert$ are from the uniform
distribution on the compact interval $\left[  0.5,3.5\right]  $ but each $\mu_{i}$ has probability $0.5$ to be negative or positive. Note that the magnitudes of the nonzero $\mu_{i}$'s are more varying than those simulated in \cite{Fan:2012}, since one major target there was to have uniformly, relatively large nonzero $\vert \mu_i \vert$'s in order to well estimate the conditional
FDR of the conditional MTP there.

The $6$ types of correlation matrix $\boldsymbol{\Sigma}=\left(\tilde{\sigma}_{ij}\right)$ are given below:
\begin{itemize}
\item ``Autoregressive'': $\tilde{\sigma}_{ij}=\rho^{\left\vert i-j\right\vert }1_{\left\{
i\neq j\right\}  }$ with $\rho=0.7$. This is the autocorrelation matrix of an
autoregressive model of order $1$. Since
$m^{-2}\left\Vert \boldsymbol{\Sigma}\right\Vert _{1}=m^{-1}\rbk{1+\frac{2\rho}%
{1-\rho}}+O\left(  m^{-2}\right)$.
\autoref{prop:SLLNMarginal} implies that we can directly implement the marginal MTP and that the SLLN will hold for the
sequences of rejections.

\item ``Block Dependence'': $\boldsymbol{\Sigma}=\operatorname{diag}\left\{  \boldsymbol{\Sigma
}^{\left(  1\right)  },\boldsymbol{\Sigma}^{\left(  2\right)  }%
,\boldsymbol{\Sigma}^{\left(  3\right)  },\boldsymbol{\Sigma}^{\left(
4\right)  }\right\}  $, where the dimension of $\boldsymbol{\Sigma}^{\left(
i\right)  }$ is $c_{i}$ for $i=1,\ldots,4$ and $c_{1}=0.1m$, $c_{2}=0.2m$, $c_{3}=0.3m$ and $c_{4}=0.4m$.
The first $2$ blocks of $\boldsymbol{\Sigma}$ are structured but the rest not, and $\boldsymbol{\Sigma}$ is a.s. singular.
The blocks are generated as follows:
$\boldsymbol{\Sigma}^{\left(  1\right)  }\left(  i,j\right)  =\rho 1_{\left\{
i\neq j\right\}  }$ with $\rho=0.7$.
$\boldsymbol{\Sigma}^{\left(  2\right)  }=\mathbf{SS}^{\top}$ with
\[
\mathbf{S}=\left(  -\frac{1}{4}\mathbf{1}_{c_{2}},\frac{1}{5}\mathbf{1}%
_{c_{2}},-\frac{1}{8}\mathbf{1}_{c_{2}},\left(  1-\sqrt{0.118125}%
\mathbf{I}_{c_{2}}\right)  \right),
\]
where $\mathbf{1}_{s}$ is a column vector of $s$ one's and $\mathbf{I}_{s}$ the $s \times s$ identity matrix.
In fact,
\[
\boldsymbol{\Sigma}^{\left(  2\right)  }=0.118125\times\mathbf{1}_{c_{2}%
}\mathbf{1}_{c_{2}}^{\top}+\left(  1-0.118125\right)  \times\mathbf{I}_{c_{1}}%
\]
since $\left(  -1/4\right)  ^{2}+\left(  1/5\right)  ^{2}+\left(  -1/8\right)
^{2}=0.118125$, and $\boldsymbol{\Sigma}^{\left(  2\right)  }$ is
the correlation matrix a Normal vector with equally correlated components and
the correlations are generated by $3$ factors.
Note that the types of dependence encoded by $\boldsymbol{\Sigma}^{\left(  1\right)  }$ and $\boldsymbol{\Sigma}^{\left(  2\right)  }$ were also used in the simulation study in \cite{Fan:2012}.
$\boldsymbol{\Sigma}^{\left(  3\right)  }$ is the sample correlation
matrix of a $c_{3}\times20$ matrix of i.i.d. Binomial random variables with
total number of trials $10$ and probability of success $0.7$. For finite $m$,
$\boldsymbol{\Sigma}^{\left(  3\right)  }$ is singular and has $20$ positive eigenvalues a.s. 
$\boldsymbol{\Sigma}^{\left(  4\right)  }$ is the sample
correlation matrix of a $c_{4}\times0.01c_{4}$ matrix of i.i.d. standard
Normal random variables. $\boldsymbol{\Sigma}^{\left(  4\right)  }$ is unstructured and a.s. singular.

\item ``Equi-correlation'': $\tilde{\sigma}_{ij} = \rho 1_{\left\{i\neq j\right\}  }$ with $\rho=0.7$.

\item ``Fractional Gaussian'':
\[
\tilde{\sigma}_{ij}=\frac{1}{2}\left[  \left(  \left\vert i-j\right\vert +1\right)
^{2H}-2\left\vert i-j\right\vert ^{2H}+\left(  \left\vert i-j\right\vert
-1\right)  ^{2H}\right]  1_{\left\{  i\neq j\right\}  }%
\]
where the Hurst index $H=0.9$. This type of long-range dependence is
pertain to fractional Gaussian noise (FGN), the increment process of
fractional Brownian motion, and has been used to model fluid
dynamics. From equation (2.13) on page 52 of \cite{Beran:1994} or identity (10) of
\cite{Zunino20086057}, we obtain%
\[
\lim_{\left\vert i-j\right\vert \rightarrow\infty}\frac{\tilde{\sigma}_{ij}}{H\left(
2H-1\right)  \left\vert i-j\right\vert ^{2H-2}}=1
\]
for $0<H<1$. So, when $H=0.9$ we have $\tilde{\sigma}_{ij}\sim0.72\left\vert i-j\right\vert ^{-0.2}$ for $i \ne j$.

\item ``Moving Average'': $\boldsymbol{\Sigma}$ is a banded matrix of bandwidth
$b=0.5m$, where $\tilde{\sigma}_{ij}=\sum_{l=1}^{b-\left\vert i-j\right\vert
}b_{l}b_{\left\vert i-j\right\vert +l}$ with $b_{i}=\frac{1}{\sqrt{b}}$ when $ 0 < \vert i-j \vert < b$.
Namely, $\boldsymbol{\Sigma}$ is the autocorrelation matrix of a moving average model of order $0.5m$.
The smallest off-diagonal nonzero entry of $\boldsymbol{\Sigma}$ is $2m^{-1}$, i.e., the weakest correlation among
two different components of the Normal random vector $\ms{\zeta}$ is $2m^{-1}$.

\item ``Unstructured Covariance'': generate an $m\times0.01m$ matrix $\mathbf{B}$ of i.i.d.
observations from the standard Normal random variable, and set
$\boldsymbol{\Sigma}$ as the sample correlation matrix of $\mathbf{B}$.
Note that $\boldsymbol{\Sigma}$ is unstructured and a.s. singular.

\end{itemize}

We briefly comment on the ranges of dependence represented by the $6$ types of
correlation matrices. When $0<\rho<1$, $0.5<H<1$ and $\left\vert
i-j\right\vert $ is large, we have
$
\rho^{\left\vert i-j\right\vert }<2m^{-1}<\left\vert i-j\right\vert ^{2H-2 }<\rho.
$
So, the ranges of dependencies, ordered from the
shortest to the longest, are roughly Autoregressive, Moving Average, Fractional
Gaussian, Equi-correlation, Block Dependence and Unstructured Covariance. However, the first $4$ types of
dependencies are all structured. Specifically, Equi-correlation and Autoregressive
are much more structured than Fractional Gaussian, and Moving Average is the least structured.
In contrast, the dependency encoded by Block Dependence or Unstructured Covariance is not structured.
The properties of the $6$ types of dependence will help compare the efficiency of how PFA and PCS
construct the major and minor vectors; see \autoref{figProjDim} for an illustration.

The simulation is implemented as follows:
\begin{enumerate}
  \item Fix a combination of $m$, $\pi_{1,m}$ and $\ms{\Sigma}$, and set the upper bound for PCS and PFA as $0.5m^{-0.4}$, i.e., the covariance matrix $\boldsymbol{\Sigma}_{\mathbf{v}_{\textrm{pfa}}}$ of the minor vector $\mb{v}_{\textrm{pfa}}$ obtained by PFA satisfies $m^{-1}\left\Vert \boldsymbol{\Sigma}_{\mathbf{v}_{\textrm{pfa}}}\right\Vert _{2} \le 0.5m^{-0.4}$,
      and the correlation matrix $\mb{R}_{\mathbf{v}_{\textrm{pcs}}}$ of the minor vector $\mb{v}_{\textrm{pcs}}$ obtained by PCS satisfies $m^{-2}\left\Vert \mb{R}_{\mathbf{v}_{\textrm{pcs}}}\right\Vert _{1} \le 0.5m^{-0.4}$. Note that the same upper bound $0.5m^{-0.4}$ is used for $m^{-1}\left\Vert \boldsymbol{\Sigma}_{\mathbf{v}_{\textrm{pfa}}}\right\Vert _{2} $ and $m^{-2}\left\Vert \mb{R}_{\mathbf{v}_{\textrm{pcs}}}\right\Vert _{1}$. Generate $\boldsymbol{\mu}$.

  \item Repeat the following $1000$ times:
   \begin{enumerate}
     \item Generate $\ms{\zeta}$ from $\mathsf{N}_{m}\left(  \boldsymbol{\mu},\boldsymbol{\Sigma}\right)$; apply PCS to obtain the major vector $\ms{\eta}_{\textrm{pcs}}$ and minor vector $\mb{v}_{\textrm{pcs}}$ such that $\ms{\zeta} = \ms{\mu}+\ms{\eta}_{\textrm{pcs}}+\mb{v}_{\textrm{pcs}}$; apply PFA to obtain the major vector $\ms{\eta}_{\textrm{pfa}}$ and minor vector $\mb{v}_{\textrm{pfa}}$  such that $\ms{\zeta} = \ms{\mu}+\ms{\eta}_{\textrm{pfa}}+\mb{v}_{\textrm{pfa}}$; the implementation of PCS is given in \autoref{secImplement} and that of PFA stated in \autoref{SecConnection}.
     \item Apply the conditional or adjusted conditional MTP to $\ms{\zeta}$ conditional on $\ms{\eta}_{\textrm{pcs}}$ and to $\ms{\zeta}$ conditional on $\ms{\eta}_{\textrm{pfa}}$ respectively. The adjusted conditional MTP is defined in \autoref{SecCondMTP}.
   \end{enumerate}
  \item Obtain the sample variance of the average number of (conditional) rejections.
  \item Repeat the previous steps to exhaust all $126$ combinations of $m$, $\pi_{1,m}$ and $\ms{\Sigma}$.
\end{enumerate}

\subsection{Implementation of PCS}
\label{secImplement}
We implement PCS through the following steps:

\begin{enumerate}
\item Obtain the spectral decomposition $\boldsymbol{\Sigma}=\sum_{i=1}%
^{m}\lambda_{i}\boldsymbol{\gamma}_{i}\boldsymbol{\gamma}_{i}^{\top}$, where
$\lambda_{i},i=1,\ldots,m$, descendingly ordered in $i$, are the eigenvalues
of $\ms{\Sigma}$ and $\boldsymbol{\gamma}_{i}$ is the eigenvector associated with
$\lambda_{i}$.

\item For each integer $k$ between $0$ and $m-1$, let $\mathbf{\tilde{Q}}%
_{k}=\sum_{i=k+1}^{m}\lambda_{i}\boldsymbol{\gamma}_{i}\boldsymbol{\gamma}%
_{i}^{\top}$ and standardize $\mathbf{\tilde{Q}}_{k}$ into a correlation matrix
$\mathbf{Q}_{k}$. Note that $\mathbf{Q}_{k}$ can be singular. Pick $\delta>0$
and a small, positive constant $C_{0}$. Find $k_{0}$, the smallest integer $k$
between $0$ and $m-1$, such that%
\begin{equation}
m^{-2}\left\Vert \mathbf{Q}_{k}\right\Vert _{1}\leq C_{0}m^{-\delta}. \label{eqb1}%
\end{equation}

\item There are three cases for $k_{0}$: (i) if $k_{0}$ does not exist, adjust
$C_{0}$ or $\delta$ so that $k_{0}$ exists; (ii) if $k_{0}=0$, set
$\boldsymbol{\eta}=\mathbf{0}$ and $\mathbf{v}=\boldsymbol{\zeta}$; in this case a conditional MTP is not needed; (iii) if $0<k_0\leq m-1$, set
$\boldsymbol{\Sigma}_{\mathbf{v}}=\sum_{i=k_0+1}^{m}\lambda_{i}%
\boldsymbol{\gamma}_{i}\boldsymbol{\gamma}_{i}^{\top}$, $\boldsymbol{\Sigma
}_{\boldsymbol{\eta}}=\sum_{i=1}^{k_0}\lambda_{i}\boldsymbol{\gamma}%
_{i}\boldsymbol{\gamma}_{i}^{\top}$, $\boldsymbol{\eta}\sim\mathsf{N}%
_{m}\left(  \mathbf{0},\boldsymbol{\Sigma}_{\boldsymbol{\eta}}\right)  $ and
$\mathbf{v}\sim\mathsf{N}_{m}\left(  \mathbf{0},\boldsymbol{\Sigma
}_{\mathbf{v}}\right)  $, so that $\boldsymbol{\zeta}=\boldsymbol{\mu
}+\boldsymbol{\eta}+\mathbf{v}$ and $\boldsymbol{\Sigma=\Sigma}%
_{\boldsymbol{\eta}}+\boldsymbol{\Sigma}_{\mathbf{v}}$.
\end{enumerate}

In the simulation study, we set $\delta=0.4$ and $C_{0}=0.5$
to demonstrate our convergence results without requiring $m$ to be extremely large since the largest value of $m$ we have is $10000$.
The same $\delta$ and $C_{0}$ values are used for PFA.
In general, the choice of $\delta$ and $C_{0}$ does
affect the performance of conditional multiple testing based on PCS, for which
smaller $C_{0}m^{2-\delta}$ makes the average number of (conditional) rejections
more concentrated around its expectation. However, $C_0$ and $\delta$ should not be considered as two tuning parameters, since in practice we will always
set a small value for $C_0 m^{-\delta}$ to upper bound $m^{-2}\left\Vert \mb{R}_{\mathbf{v}_{\textrm{pcs}}}\right\Vert _{1}$ or $m^{-2}\left\Vert \mathbf{Q}_{k}\right\Vert _{1}$ to implement PCS, as was done to upper bound
$m^{-1}\left\Vert \boldsymbol{\Sigma}_{\mathbf{v}_{\textrm{pfa}}}\right\Vert _{2} $ to implement PFA.

Since, to determine $k_0$ in Step 2 of the implementation of PCS, computing the correlation matrix $\mathbf{Q}_{k}$ for each $k=0,\ldots,m$ can take some computational time when $m$ is very large, in the simulation, we let $k$ run through a nonlinear sequence of distinct numbers $1,2,\ldots,10, \tilde{k}_1, \ldots, \tilde{k}_{l}, m $, where each $\tilde{k}_i$ is the integer part of $10\times q^i$ for $i=1, \ldots,l$, $q$ is a prespecified positive real number that controls the length of the sequence, and $l$ is the integer part of $\log_{q}\rbk{10^{-1}m}$. Specifically, $q=1.01$ is set when $m <10^3$, and $q=1.02$ when $10^3 \le m \le 10^4$. The first $10$ consecutive numbers in the sequence is to cover types of dependence where only a very smaller number of eigenvectors of the covariance matrix of the Normal random vector is needed to construct the major vector to achieve PCS. Better performance of PCS in terms of more concentrated average number of conditional rejections may be obtained if we implement exactly Step 2.

\subsection{Summary of simulation results}\label{simRes}

\begin{figure}[t!]
\centering
\includegraphics[width=.9\textwidth]{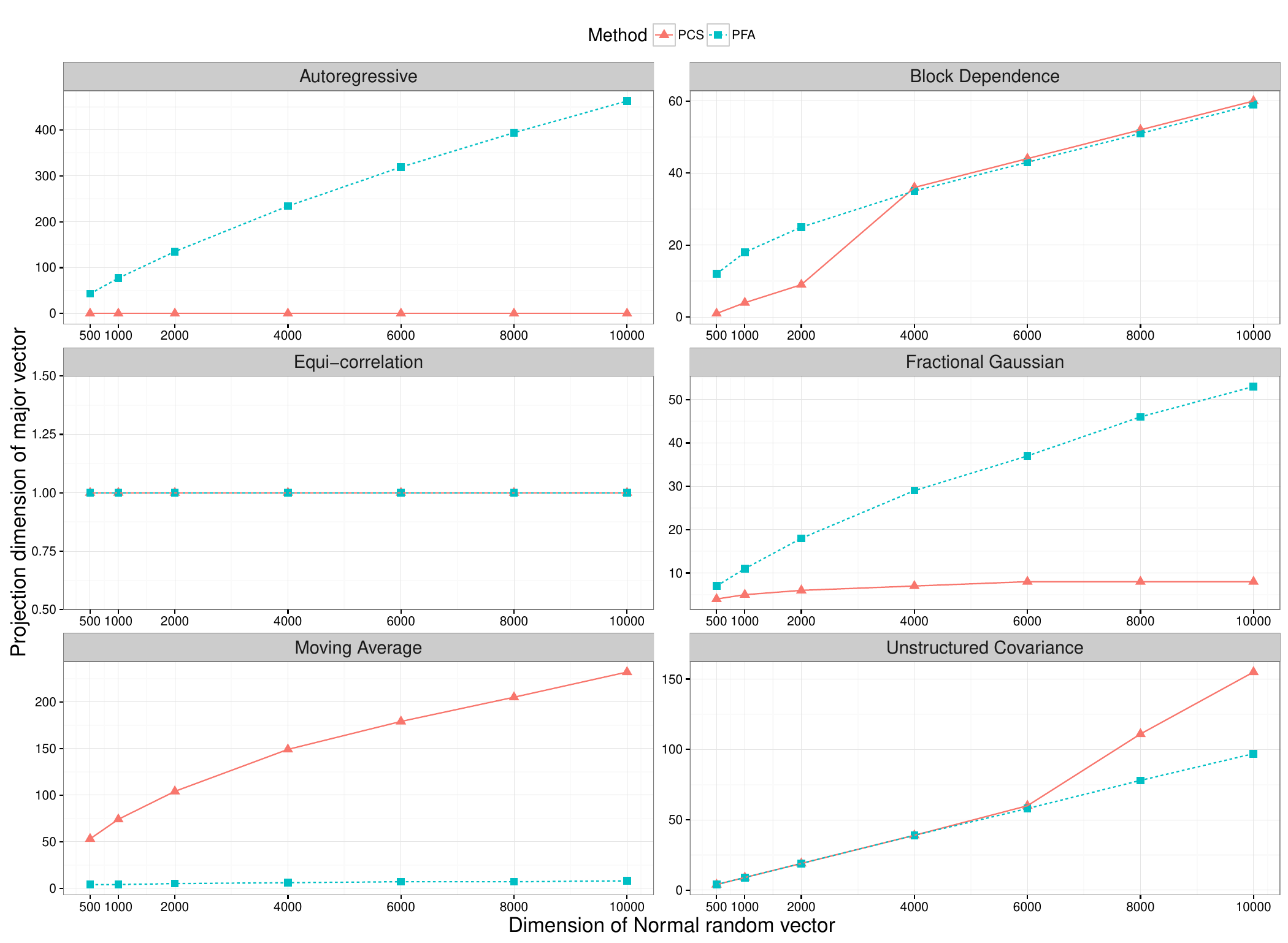}
\caption[]{The projection dimension $\chi\rbk{\ms{\eta}}$ of the major vector $\ms{\eta}$ for different types of dependence as the dimension of the Normal random vector $\ms{\zeta}$ changes. $\chi\rbk{\ms{\eta}}$  is the number of eigenvectors of the covariance matrix $\ms{\Sigma}$ of $\ms{\zeta}$ used to construct $\ms{\eta}$. PCS (denoted by triangle) adaptively determines $\chi\rbk{\ms{\eta}}$ according to the complexity of the dependence structure $\ms{\Sigma}$, i.e., a more structured dependence leads to a smaller projection dimension. In contrast, the projection dimension based on PFA (denoted by square) may be excessive (see Fractional Gaussian and Autoregressive) or insufficient (see Moving Average).}
\label{figProjDim}
\end{figure}

We first assess the efficiency of PCS and PFA in constructing the major vector $\ms{\eta}$ in terms of ``projection dimension''. The projection dimension $\chi\rbk{\ms{\eta}}$ of $\ms{\eta}$ is the number of eigenvectors of the covariance matrix $\ms{\Sigma}$ of $\ms{\zeta}$ used to construct $\ms{\eta}$.
$\chi\rbk{\ms{\eta}}$ for PCS and for PFA are given in \autoref{figProjDim}.
PCS adaptively determines $\chi\rbk{\ms{\eta}}$ according to the complexity of the dependence structure encoded by $\ms{\Sigma}$, i.e., a more structured dependency leads to a smaller projection dimension.
For example, for Autoregressive, $\ms{\Sigma}$ already has principal correlation structure. So, the conditional MTP becomes the marginal MTP and $\chi\rbk{\ms{\eta}}=0$. Further, among the other $3$ types of structured dependence, Equi-correlation is more structured than Fractional Gaussian, and Moving Average is the least structured. Accordingly, $\chi\rbk{\ms{\eta}}$ based on PCS is the largest for Moving Average, smaller for Fractional Gaussian, and the smallest for Equi-correlation. In contrast, $\chi\rbk{\ms{\eta}}$ based on PCS for unstructured dependence such as Block Dependence and Unstructured Covariance is usually larger than those for the structured dependencies.
However, $\chi\rbk{\ms{\eta}}$ based on PFA does not seem to adapt to the complexity of dependence. For example, it may be excessive (see Fractional Gaussian and Autoregressive) or insufficient (see Moving Average).

\begin{figure}[t!]
\centering
\includegraphics[width=.9\textwidth]{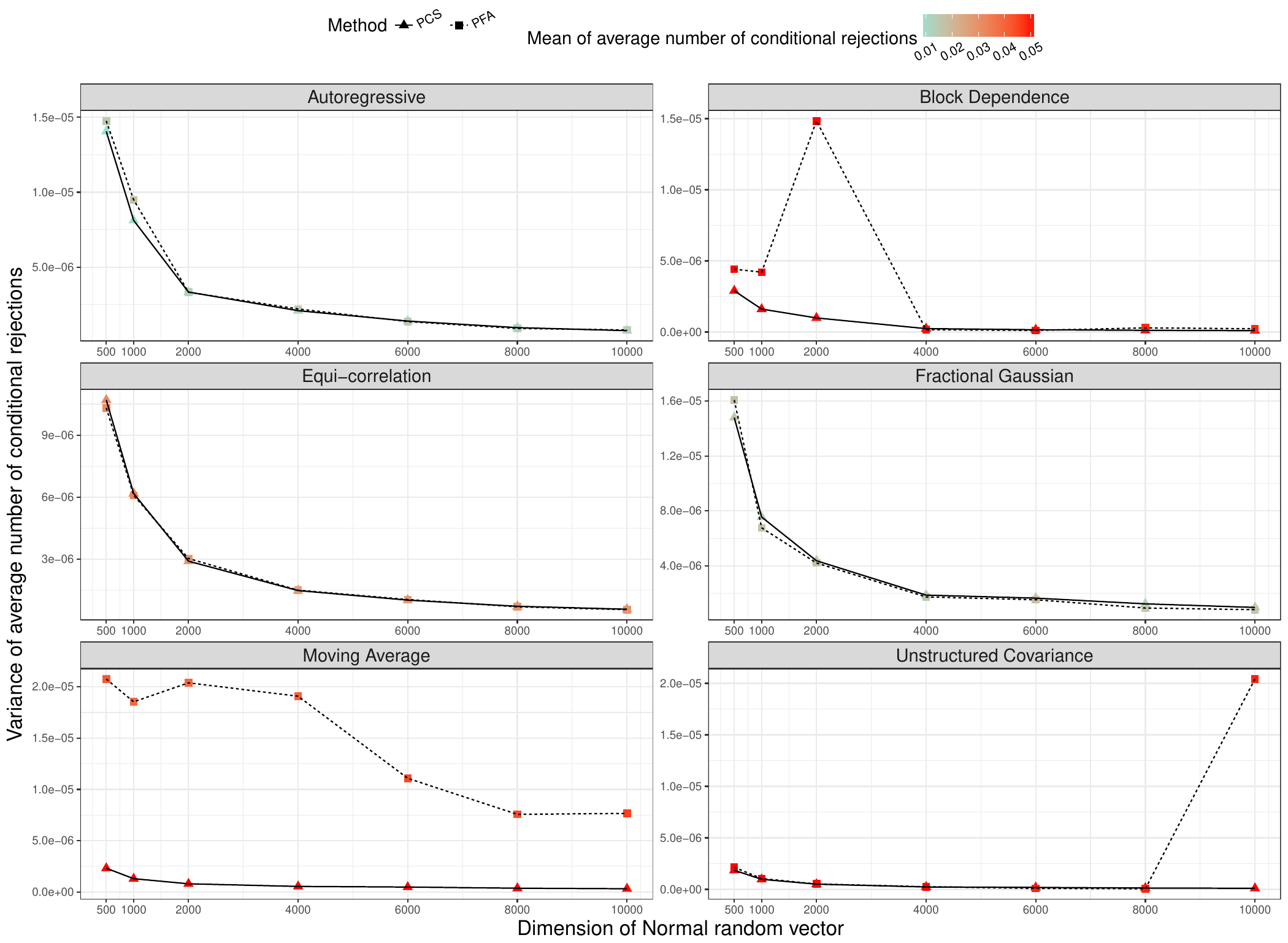}
\caption[]{Variance of the average number of conditional rejections $m^{-1}R_m\rbk{t\vert\ms{\eta}}$ of the adjusted conditional MTP in the dense regime. For PCS (denoted by triangle) we see a steady trend that the variances converge to $0$ as $m$ increase. However, for PFA (denoted by square), as $m$ increases, the variances do not necessarily show a trend of converging to $0$ (see Block Dependence or Unstructured Covariance), or they can show a trend of slow convergence to $0$ (see Moving Average).
}
\label{figAvgRej}
\end{figure}

\begin{figure}[t!]
\centering
\includegraphics[width=1\textwidth]{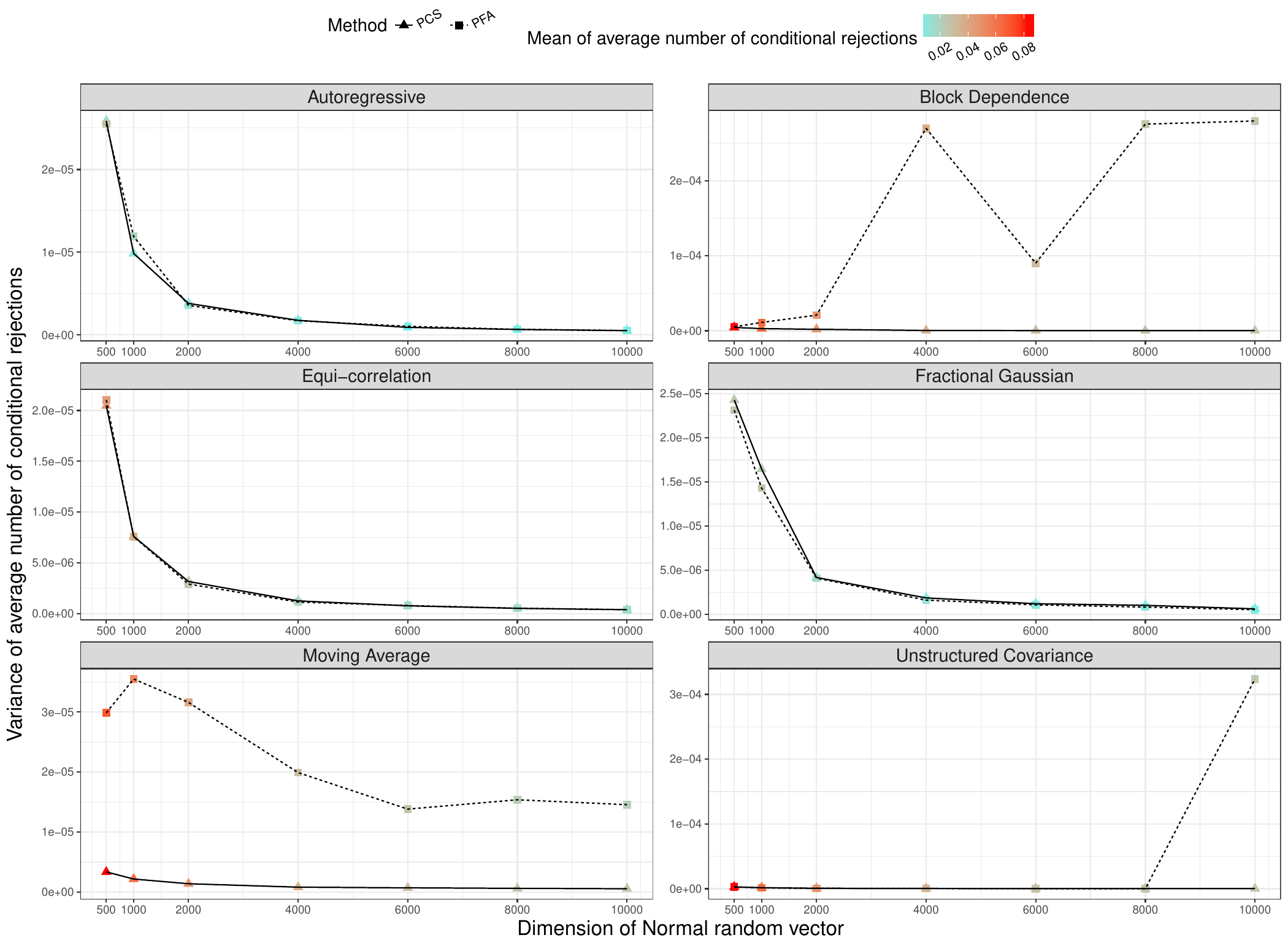}
\caption[]{Variance of the average number of conditional rejections $m^{-1}R_m\rbk{t\vert\ms{\eta}}$ of the adjusted conditional MTP in the moderately sparse regime. For PCS (denoted by triangle) we see a steady trend that the variances converge to $0$ as $m$ increase. However, for PFA (denoted by square), as $m$ increases, the variances do not necessarily show a clear trend of converging to $0$ (see Block Dependence or Unstructured Covariance).}
\label{figAvgRejms}
\end{figure}

\begin{table}[t!]
\begin{center}
\begin{tabular}{l|l|l|l|l}
  \hline
$m$ & Dependence  & $\vartheta_0$ & $\lambda_0$ & $\kappa_0$ \\
  \hline
  2000 & Block Dependence &         4 &  800 & 0.00000100 \\
   4000 & Block Dependence &   1438800 & 1600 & 0.08992500 \\
   6000 & Block Dependence &   3238200 & 2400 & 0.08995000 \\
   6000 & Unstructured Covariance &  35982002 &    0 & 0.99950006 \\
   8000 & Block Dependence &  10236800 & 2400 & 0.15995000 \\
   8000 & Unstructured Covariance &  63976000 &    1 & 0.99962500 \\
  10000 & Block Dependence &  15988002 & 3000 & 0.15988002 \\
  10000 & Unstructured Covariance &    283908 &    0 & 0.00283908 \\
   \hline
\end{tabular}
\caption{Minor vector $\mb{v}=\rbk{v_1,\cdots,v_m}^{\top}$ obtained by PFA: the number $\vartheta_0$ of distinct pairs $\rbk{v_i,v_j}$ whose absolute correlations (``absolute correlation'' is the absolute value of correlation) are at least $0.99999 = 1- 10^{-5}$, the number $\lambda_0$ of components of $\mb{v}$ whose standard deviations are no larger than $10^{-5}$, and the proportion $\kappa_0 = \vartheta_0/m^2$. This table is for \autoref{figAvgRejms}.
}
\label{TbLinear}
\end{center}

\end{table}

\begin{figure}[t!]
\centering
\includegraphics[width=1\textwidth]{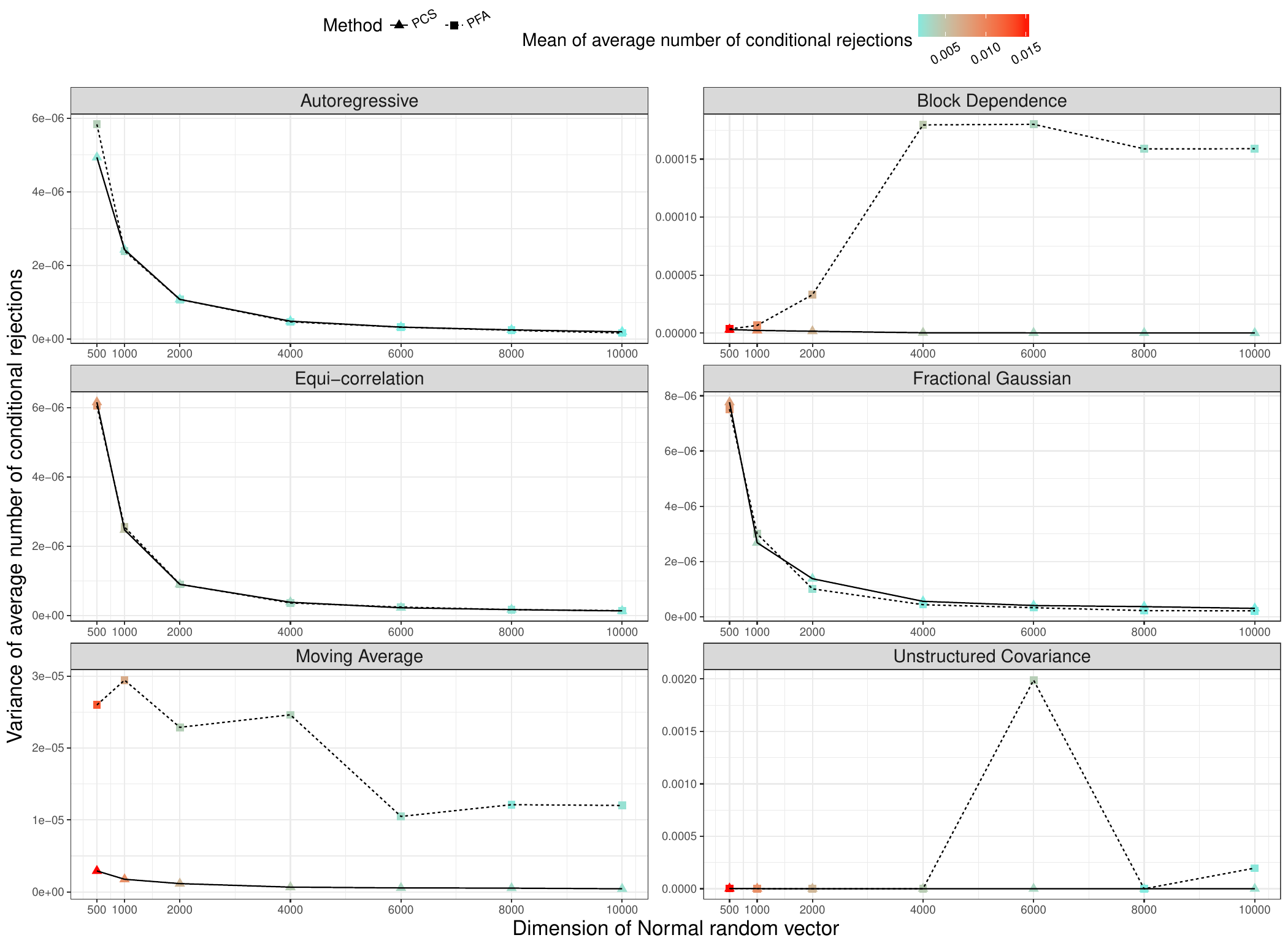}
\caption[]{Variance of the average number of conditional rejections $m^{-1}R_m\rbk{t\vert\ms{\eta}}$ of the adjusted conditional MTP in the very sparse regime. For PCS (denoted by triangle) we see a steady trend that the variances converge to $0$ as $m$ increase. However, for PFA (denoted by square), as $m$ increases, the variances do not necessarily show a clear trend of converging to $0$ (see Block Dependence or Unstructured Covariance).}
\label{figAvgRejvs}
\end{figure}

Secondly, we compare convergence results of $m^{-1}R_m\rbk{t\vert\ms{\eta}}, m \ge 1$ based on PCS to those based on PFA.
For each simulated type of dependence, the sequence of variances of $m^{-1}R_m\rbk{t\vert\ms{\eta}}$ based on PCS shows a strong trend of convergence to $0$ as $m$ increases, whereas that based on PFA does not necessarily; see \autoref{figAvgRej}, \autoref{figAvgRejms} and \autoref{figAvgRejvs} for the adjusted conditional MTP.
For PFA, the variance of $m^{-1}R_m\rbk{t\vert\ms{\eta}}$ under Block Dependence or Unstructured Covariance can be unstable and inflated. A major reason for this is the large number $\vartheta_0$ of almost linearly dependent pairs of components of the minor vector or a large number $\lambda_0$ of components of the minor vector whose variances are very small, both produced by PFA; see \autoref{TbLinear}. We have pointed out in \cite{Chen:2014SLLN} that the variance of $m^{-1}R_m\rbk{t\vert\ms{\eta}}$ based on PFA may be inflated when $\vartheta_0$ or $\lambda_0$ is large, and our simulation provides evidence on this.
Moreover, the variance of $m^{-1}R_m\rbk{t\vert\ms{\eta}}$ based on PFA may display a trend of very slow convergence to $0$ as $m$ increase; see, e.g., Moving Average for the adjust conditional MTP. This is mainly because PFA insufficiently determines the projection dimension of the major vector as shown in \autoref{figProjDim}. Lastly, the variances based on PFA and PCS for Autoregressive or Equi-correlation have the same profile, because for Autoregressive PCS uses more eigenvectors to construct the major vector and for Equi-correlation both PCS and PFA use one eigenvector to construct the major vector as shown in \autoref{figProjDim}.

For the conditional MTP, the variances of the sequence $m^{-1}R_m\rbk{t\vert\ms{\eta}}, m \ge 1$ based on PCS and PFA are very similar. So, we omit presenting them here.

\section{Discussion}

\label{SecDiscussion}

For the Normal means problem under dependence with a known covariance matrix, we have classified via the concept of ``principal correlation structure (PCS)'' and the $l_1$-norm of the correlation matrix of the Normal random variables the type of dependence under which the SLLN holds for the numbers of rejections and false rejections and the FDP of an MTP. Further, our simulation study suggests that multiple testing Normal means based on PCS is more stable than that based on the PFA of \cite{Fan:2012}.
The validity of our work, excluding the implementation of PCS via the spectral decomposition of the covariance matrix of the Normal random variables, only requires these random variables or the entries of the minor vector in the additive decomposition (\ref{eqModel}) to be bivariate Normal.

When the covariance matrix of the Normal random variables is unknown but can be consistently estimated, the SLLN we have established for an MTP can be easily downgraded to the weak law of large numbers (WLLN). However, estimating a covariance matrix is not the focus and is beyond the scope of the current work. Since for the SLLN associated the marginal MTP, the sufficient condition we have provided is $m^{-2}\left\Vert \mathbf{R}\right\Vert _{1}=O\left(  m^{-\delta
}\right)$ for some $\delta>0$, it is necessary to develop tests to check the order of $\left\Vert \mathbf{R}\right\Vert _{1}$.
Finally, it is possible to employ the universal comparison result provided by \autoref{ThmSLLN} and Laguerre polynomials associated with Chi-squared distributions to study the SLLN for multiple testing location parameters of Student t-statistics. We leave these two tasks to future work.

\appendix

\section{Proofs}\label{secProofs}
Let $C>0$ be a generic constant that can assume different (and appropriate) values at different occurrences.
In this section, we provide the proofs of \autoref{BndHermitePoly}, \autoref{ThmSLLN}, \autoref{prop:SLLNMarginal}, \autoref{ThmSLLNAMarg} and \autoref{ThmSLLNStrongest}.

\subsection{Proof of \autoref{BndHermitePoly}}
\label{ProofBndHermite}

Recall the $n$th Hermite polynomial defined by \cite{Mehler:1866} and used in \autoref{sec:Mehler} as%
\begin{equation*}
H_{n}\left(  x\right)  =\left(  -1\right)  ^{n}\frac{1}{\phi\left(  x\right)
}\frac{d^{n}}{dx^{n}}\phi\left(  x\right)  
\end{equation*}
where%
\[
\phi\left(  x\right)  =\left(  2\pi\right)  ^{-1/2}\exp\left(  -x^{2}%
/2\right)  .
\]
We aim to prove \eqref{eqBoundHermite}, i.e., for some constant $K_{0}>0$ independent of $n$,
\begin{equation*}
\left\vert e^{-x^{2}/2}H_{n}\left(  x\right)  \right\vert \leq K_{0}\sqrt
{n!}n^{-1/12}e^{-x^{2}/4}\text{ \ for any\ }x\in\mathbb{R}.
\end{equation*}

In order to show this, we
need to use contents from three sections of \cite{Szego:1939}: Section 1.81 on
the Airy function, Section 1.71 on Bessel functions, and Section 8.91 on the
asymptotic properties of Laguerre and Hermite polynomials defined there slightly differently. We describe the strategy of proof first,
which has three consecutive parts:

\textbf{Part 1}: Show that the Airy function $A\left(  x\right)  $ defined in
Section 1.81 and appearing in identity (8.91.10) on page 236 of
\cite{Szego:1939} is uniformly bounded, i.e.,%
\begin{equation}
\sup_{x\in\mathbb{R}}\left\vert A\left(  x\right)  \right\vert <\infty.
\label{eq13}%
\end{equation}

\textbf{Part 2}: The $n$th (physicists') Hermite polynomial $\hat{H}_{n}\left(  x\right)  $
is defined by (5.5.3) on page 102 of \cite{Szego:1939} as%
\begin{equation*}
e^{-x^{2}}\hat{H}_{n}\left(  x\right)  =\left(  -1\right)  ^{n}\left(
\frac{d}{dx}\right)  ^{n}e^{-x^{2}}\text{.} 
\end{equation*}
We will show the relationship%
\begin{equation}
\hat{H}_{n}\left(  x\right)  =2^{n/2}H_{n}\left(  \sqrt{2}x\right)
\label{eq9}%
\end{equation}
for any $x\in\mathbb{R}$.

\textbf{Part 3}: The identity (8.91.10) on page 236 of \cite{Szego:1939} on
$\hat{H}_{n}\left(  x\right)  $ reads%
\begin{equation}
\max_{x\in\mathbb{R}}e^{-x^{2}/2}\left\vert \hat{H}_{n}\left(  x\right)
\right\vert \cong\left(  2^{n}n!\right)  ^{1/2}2^{1/4}3^{1/2}\pi
^{-3/4}n^{-1/12}\max_{t\in\mathbb{R}}A\left(  t\right)  , \label{eq3}%
\end{equation}
where the notation $\cong$, defined in paragraph 6 of page 1 of
\cite{Szego:1939}, means that the ratio between the two sequences on both
sides of $\cong$ converges to $1$. Since we have
already shown that the Airy function $A\left(  x\right)  $ is uniformly
bounded, we see from (\ref{eq13}) and (\ref{eq3}) that%
\begin{equation}
\max_{x\in\mathbb{R}}e^{-x^{2}/2}\left\vert \hat{H}_{n}\left(  x\right)
\right\vert \leq K\left(  2^{n}n!\right)  ^{1/2}n^{-1/12} \label{eq6}%
\end{equation}
for some finite, positive constant $K$ independent of $n$. Clearly, (\ref{eq6}) implies
inequality (30) of \cite{Hille1926}, i.e.,%
\begin{equation}
\left\vert \hat{H}_{n}\left(  x\right)  \right\vert \leq K_{0}2^{n/2}\sqrt
{n!}n^{-1/12}e^{x^{2}/2} \label{eq5}%
\end{equation}
holds for some constant $K_{0}>0$ independent of $n$.

Finally, plugging (\ref{eq9}) into (\ref{eq5}) gives%
\[
\left\vert 2^{n/2}H_{n}\left(  \sqrt{2}x\right)  \right\vert \leq
K2^{n/2}\sqrt{n!}n^{-1/12}e^{x^{2}/2},
\]
which is equivalent to%
\begin{equation}
\left\vert H_{n}\left(  x\right)  \right\vert \leq K\sqrt{n!}n^{-1/12}%
e^{x^{2}/4}. \label{eq8}%
\end{equation}
Now multiple both sides of (\ref{eq8}) by $e^{-x^{2}/2}$, we get exactly
inequality \eqref{eqBoundHermite} in the main text.

Now we provide details for {\bf Parts 1} and {\bf 2}.

\textbf{Step 1}: Let us show that the Airy function $A\left(  x\right)  $ defined
in Section 1.81 (pages 18 and 19) of \cite{Szego:1939} is uniformly bounded in
$x\in\mathbb{R}$. The facts we will use to show this are contained in Section
1.81 and Section 1.71 of \cite{Szego:1939}. So, when we state them we will not
explicitly mention the source \cite{Szego:1939} every time. Modulo a constant factor, the Airy function
$A\left(  x\right)  $ is defined by (1.81.1) and (1.81.4), i.e.,
\[
A\left(  x\right)  =k\left(  x\right)  +l\left(  x\right)  ,
\]
where%
\[
k\left(  x\right)  =\frac{\pi}{3}\left(  x/3\right)  ^{1/2}J_{-1/3}\left(
2\left(  x/3\right)  ^{3/2}\right)
\]
and%
\[
l\left(  x\right)  =\frac{\pi}{3}\left(  x/3\right)  ^{1/2}J_{1/3}\left(
2\left(  x/3\right)  ^{3/2}\right)  ,
\]
where $J_{\cdot}\left(  x\right)  $ is the Bessel function of the first kind.
From (1.81.5), i.e., %
\[
A\left(  x\right)  \cong2^{-1}3^{-1/4}\pi^{1/2}\left\vert x\right\vert
^{1/2}\exp\left\{  -2\left(  \left\vert x\right\vert /3\right)  ^{1/2}%
\right\}  \text{ \ as \ }x\rightarrow-\infty,
\]
we see $\lim_{x\rightarrow-\infty}A\left(  t\right)  =0$.

By the identity (1.71.11) on page 16
in Section 1.71 for the asymptotic order of Bessel functions of the first kind
as $x\rightarrow+\infty$, i.e.,%
\[
J_{\alpha}\left(  z\right)  =O\left(  z^{-1/2}\right)  \text{ \ as
\ }z\rightarrow+\infty
\]
for any real $\alpha$, we see that%
\[
\left\vert A\left(  x\right)  \right\vert \leq C\pi x^{1/2}\left(
x^{3/2}\right)  ^{-1/2}=C\pi x^{-1/4}%
\]
as $x\rightarrow+\infty$. Therefore, $\lim_{x\rightarrow+\infty}\left\vert
A\left(  t\right)  \right\vert =0$. However, we already have $\lim
_{x\rightarrow-\infty}A\left(  t\right)  =0$ and $A\left(  x\right)  $ is
continuously differentiable in $x$ for all $x\in\mathbb{R}$. Thus, (\ref{eq13})
holds, i.e., $A\left(  x\right)  $ is uniformly bounded in $x$ for
$x\in\mathbb{R}$.

\textbf{Step 2}: Recall
\[
\hat{H}_{n}\left(  x\right)  =\left(  -1\right)  ^{n}e^{-x^{2}}\frac{d^{n}%
}{dx^{n}}e^{-x^{2}}%
\]
and%
\[
H_{n}\left(  x\right)  =\left(  -1\right)  ^{n}e^{-x^{2}/2}\frac{d^{n}}%
{dx^{n}}e^{-x^{2}/2}.
\]
Now we will show (\ref{eq9}), i.e.,%
\[
\hat{H}_{n}\left(  x\right)  =2^{n/2}H_{n}\left(  \sqrt{2}x\right)  \text{
\ for }n=0,1,....
\]
Obviously, $\hat{H}_{0}\left(  x\right)  =H_{0}\left(  \sqrt{2}x\right)  $.
Let $x=\frac{\sqrt{2}}{2}\tilde{x}$ for $\tilde{x}\in\mathbb{R}$. Then
$e^{-x^{2}}=e^{-\tilde{x}^{2}/2}$. Once we justify%
\begin{equation}
\frac{d^{n}}{dx^{n}}e^{-\tilde{x}^{2}/2}=2^{n/2}\frac{d^{n}}{d\tilde{x}^{n}%
}e^{-\tilde{x}^{2}}\text{ \ for }n=1,2,.... \label{eq10}%
\end{equation}
Then%
\begin{align}
\hat{H}_{n}\left(  \frac{\sqrt{2}}{2}\tilde{x}\right)   &  =\left(  -1\right)
^{n}e^{-\tilde{x}^{2}/2}\frac{d^{n}}{dx^{n}}\exp\left(  -\tilde{x}%
^{2}/2\right)\nonumber\\
&  =2^{n/2}\left(  -1\right)  ^{n}e^{-\tilde{x}^{2}/2}\frac{d^{n}}{d\tilde
{x}^{n}}e^{-\tilde{x}^{2}/2}
  =2^{n/2}H_{n}\left(  \tilde{x}\right) \nonumber , 
\end{align}
which is equivalent to (\ref{eq9}), i.e., $\hat{H}_{n}\left(  x\right)
=2^{n/2}H_{n}\left(  \sqrt{2}x\right)  $.

To show (\ref{eq10}), we use induction. Clearly, (\ref{eq10}) holds
automatically for $n=0$. For $n=1$, we have by chain rule%
\[
\frac{d}{dx}e^{-\tilde{x}^{2}/2}=\frac{d\tilde{x}}{dx}\times\frac{d}%
{d\tilde{x}}e^{-\tilde{x}^{2}}=2^{1/2}\frac{d}{d\tilde{x}}e^{-\tilde{x}^{2}}.
\]
Now suppose%
\[
\frac{d^{n}}{dx^{n}}e^{-\tilde{x}^{2}/2}=2^{n/2}\frac{d^{n}}{d\tilde{x}^{n}%
}e^{-\tilde{x}^{2}}%
\]
as the induction hypothesis. Then%
\begin{align*}
\frac{d^{n+1}}{dx^{n+1}}e^{-\tilde{x}^{2}/2}  &  =\frac{d}{dx}\left(
\frac{d^{n}}{dx^{n}}e^{-\tilde{x}^{2}/2}\right)  =\frac{d}{dx}\left(
2^{n/2}\frac{d^{n}}{d\tilde{x}^{n}}e^{-\tilde{x}^{2}}\right) \\
&  =2^{n/2}\frac{d}{dx}\frac{d^{n}}{d\tilde{x}^{n}}e^{-\tilde{x}^{2}}
  =2^{n/2}\frac{d\tilde{x}}{dx}\times\frac{d}{d\tilde{x}}\frac{d^{n}}%
{d\tilde{x}^{n}}e^{-\tilde{x}^{2}}
  =2^{\left(  n+1\right)  /2}\frac{d^{n+1}}{d\tilde{x}^{n+1}}e^{-\tilde
{x}^{2}}%
\end{align*}
Therefore, (\ref{eq10}) holds, and so does (\ref{eq9}).
This completes the whole proof.

\subsection{Proof of \autoref{ThmSLLN}}
\label{secBndVar}

Let $\xi_{ij}=\mathrm{cov}_{\mathbf{v}}\left(  X_{i},X_{j}\right)  $ and
$\boldsymbol{\Sigma}_{\mathbf{v}}=\left(  q_{ij}\right)  _{m\times m}$ be the
covariance matrix of $\mathbf{v}$. Note that any $v_{i}$ whose standard deviation
$\sigma_{i,m}=0$ contributes nothing to $\mathbb{V}_{\mathbf{v}}\left[
m^{-1}R_{m}\left(  t|\boldsymbol{\eta}\right)  \right]  $, so we only need to
deal with $v_{i}$ whose $\sigma_{i,m}>0$. First, we deal
with linearly dependent pairs $\left(  v_{i},v_{j}\right)  $ with $i\neq j$,
i.e., pairs $\left(  i,j\right)  \in$ $E_{2,m}$, where $E_{2,m}$ is defined in
(\ref{eqSets}). Clearly, $\left\vert \xi_{ij} \right\vert \le C \left\vert \rho_{ij}\right\vert$
for each $\left(i,j\right) \in E_{2,m}$, where $\rho_{ij}$ is the correlation between $v_i$ and $v_j$. Since $\left\vert E_{2,m}\right\vert =O\left(  m^{2-\delta
}\right)  $, we have
\[
m^{-2}\sum\nolimits_{\left(  i,j\right)  \ \in E_{2,m}}\left\vert \xi
_{ij}\right\vert \leq Cm^{-\delta}.
\]
Further, $m^{-2}\sum_{i=1}^{m}\mathbb{V}_{\mathbf{v}}\left[  X_{i}\right]
=O\left(  m^{-1}\right)  $. So, $I_{1}$ defined in (\ref{eqVarPa}) satisfies
$\left\vert I_{1}\right\vert =O\left(  m^{-\min\left\{  \delta,1\right\}
}\right)  $.

Next, we consider pairs $\left(  v_{i},v_{j}\right)  $ with $i\neq j$ that are
not linearly dependent, i.e., pairs $\left(  i,j\right)  \in E_{1,m}$, where
$E_{1,m}$ is defined in (\ref{eqSets}). Recall $c_{1,i}=\sigma_{i,m}%
^{-1}r_{1,i}$ and let $\Psi_{m}=\sum_{\left(  i,j\right)  \in E_{1,m}}%
\xi_{ij}$. For the rest of the proof, we focus on the case of one-sided
p-values since the case of two-sided ones can be dealt with similarly.

\textbf{Case 1}: one-sided p-values. Then \autoref{LmVarianceNormal} and
(\ref{eqBoundHermite}) imply
\begin{equation}\label{BndCondCov}
\left\vert \xi_{ij}\right\vert \leq \frac{\left\vert q_{ij}\right\vert }{\sigma
_{i,m}\sigma_{j,m}}\sum_{n=1}^{\infty}n^{-7/6}\left\vert \rho_{ij}\right\vert
^{n-1}\exp\left(  -4^{-1}c_{1,i}^{2}\right)  \exp\left(  -4^{-1}c_{1,j}%
^{2}\right) \le C \left\vert \rho_{ij}\right\vert
\end{equation}
(which is the universal comparison inequality (\ref{IneqComp})) and
\[
\left\vert \Psi_{m}\right\vert \leq\tilde{\Psi}_{m}=m^{-2}\sum_{\left(
i,j\right)  \in E_{1,m}}\frac{\left\vert q_{ij}\right\vert }{\sigma
_{i,m}\sigma_{j,m}}\sum_{n=1}^{\infty}n^{-7/6}\left\vert \rho_{ij}\right\vert
^{n-1}\exp\left(  -4^{-1}c_{1,i}^{2}\right)  \exp\left(  -4^{-1}c_{1,j}%
^{2}\right).
\]
So,%
\begin{equation*}
\tilde{\Psi}_{m}\leq Cm^{-2}\sum_{\left(  i,j\right)  \in E_{1,m}}%
\frac{\left\vert q_{ij}\right\vert }{\sigma_{i,m}\sigma_{j,m}}\exp\left(
-4^{-1}c_{1,i}^{2}\right)  \exp\left(  -4^{-1}c_{1,j}^{2}\right)  .
\end{equation*}
If $\sigma_{0}>0$, then $\left\vert E_{0}\right\vert =\varnothing$ and%
\begin{equation}
\tilde{\Psi}_{m}\leq Cm^{-2}\sum_{\left(  i,j\right)  \in E_{1,m}}\left\vert
q_{ij}\right\vert \leq m^{-2}\left\Vert \boldsymbol{\Sigma}_{\mathbf{v}%
}\right\Vert _{1}=O\left(  m^{-\delta}\right)  \label{eqBndOriginal}%
\end{equation}
by the assumption, where the upper bound in (\ref{eqBndOriginal}) is
independent of $\boldsymbol{\eta}$. This justifies (\ref{eqKeyVarBound}). If
$\sigma_{0}=0$, then $\left\vert E_{0}\right\vert \neq\varnothing$. Recall
$r_{1,i}=\tilde{t}-\mu_{i}-\eta_{i}$, $r_{2,i}=-\infty$, and
$G_{m,\boldsymbol{\eta}}\left(  t,\varepsilon_{m}\right)  $ in
(\ref{eqExceptionB}), i.e.,%
\[
G_{m,\boldsymbol{\eta}}\left(  t,\varepsilon_{m}\right)  =\bigcup
\nolimits_{i\in E_{0,m}}\left\{  \omega \in\Omega:\min\left\{
\left\vert r_{1,i}\right\vert ,\left\vert r_{2,i}\right\vert \right\}
<\varepsilon_{m}\right\}  .
\]
Using the fact that
\begin{equation*}
\max_{x>0}xe^{-4^{-1}x^{2}y^{2}}=\sqrt{2}y^{-1}e^{-1/2}\text{ for any }y>0,
\end{equation*}
we obtain, on the complement $D_{m,\boldsymbol{\eta}}\left(  t,\varepsilon
_{m}\right)  $ of $G_{m,\boldsymbol{\eta}}\left(  t,\varepsilon_{m}\right)
$,
\[
\sigma_{i,m}^{-1}\exp\left(  -4^{-1}c_{l,i}^{2}\right)  \leq2e^{-1/2}%
\left\vert r_{l,i}\right\vert ^{-1}\leq2\varepsilon_{m}^{-2}%
\]
and%
\[
\tilde{\Psi}_{m}\leq2e^{-1}m^{-2}\sum_{\left(  i,j\right)  \in E_{1,m}%
}\left\vert q_{ij}\right\vert \left\vert r_{1,i}\right\vert ^{-1}\left\vert
r_{1,j}\right\vert ^{-1}\sum_{n=1}^{\infty}n^{-7/6}\left\vert \rho
_{ij}\right\vert ^{n-1}.
\]
This implies $\left\vert \Psi_{m}\right\vert \leq\tilde{\Psi}_{m}\leq Cm^{-\delta }\varepsilon_{m}^{-2}$
and (\ref{eqKeyVarBoundA}) for each $\ms{\eta}\rbk{\omega}$ with $ \omega \in D_{m,\boldsymbol{\eta}}\left(t,\varepsilon_{m}\right)$.

\textbf{Case 2}: two-sided p-values. In this case, $d_{n}\left(
c,c^{\prime\prime}\right)  $ defined in \autoref{LmVarianceNormal} satisfies
\[
\left\vert d_{n}\left(  c,c^{\prime\prime}\right)  \right\vert \leq\left\vert
H_{n}\left(  c\right)  \phi\left(  c\right)  \right\vert +\left\vert
H_{n}\left(  c^{\prime\prime}\right)  \phi\left(  c^{\prime\prime}\right)
\right\vert ,
\]
and the arguments for \textbf{Case 1} lead to \eqref{BndCondCov} and the same conclusions on
$\mathbb{V}_{\mathbf{v}}\left[  m^{-1}R_{m}\left(  t|\boldsymbol{\eta
}   \right)  \right]  $. This completes the proof.

\subsection{Proof of \autoref{prop:SLLNMarginal}} \label{ProofSLLNMarginal}

The main result we rely on is quoted as follows:

\begin{lemma}
[Theorem 1 of \cite{Lyons:1988}]\label{LemmaLyons}Let $\left\{  \chi_{n}\right\}
_{n=1}^{\infty}$ be a sequence of complex-valued random variables such
that $\mathbb{E}\left[  \left\vert \chi_{n}\right\vert ^{2}\right]  \leq1$.
Set \thinspace$Q_{N}=N^{-1}\sum\nolimits_{n=1}^{N}\chi_{n}$. If \ $\left\vert
\chi_{n}\right\vert \leq1$ a.s. and%
\begin{equation}
\sum\nolimits_{N=1}^{\infty}N^{-1}\mathbb{E}\left[  \left\vert Q_{N}%
\right\vert ^{2}\right]  <\infty, \label{eqCondLyons}%
\end{equation}
then $\lim_{N\rightarrow\infty}Q_{N}=0$ a.s.
\end{lemma}
\noindent A sufficient condition for the SLLN to hold for $\left\{  \chi_{n}\right\}
_{n=1}^{\infty}$ is $\mathbb{E}\left[  \left\vert Q_{m}\right\vert
^{2}\right]  =O\left(  m^{-\delta}\right)  $ for some $\delta>0$, which
implies (\ref{eqCondLyons}).

Now we state the arguments.
When $\boldsymbol{\zeta}=\left(  \zeta_{1},\ldots,\zeta_{m}\right)  ^{\top}%
\sim\textrm{N}_{m}\left(  \boldsymbol{\mu},\boldsymbol{\Sigma}\right)  $ with
$\boldsymbol{\Sigma}=\left(  \tilde{\sigma}_{ij}\right)  $ and its correlation
matrix $\mathbf{R}$ satisfies $m^{-2}\left\Vert \mathbf{R}\right\Vert
_{1}=O\left(  m^{-\delta}\right)  $ for some $\delta>0$, we can directly
implement the marginal MTP. Recall the one-sided p-value $p_{i}=1-\Phi\left(
\tilde{\sigma}_{ii}^{-1/2}\left\vert \zeta_{i}\right\vert \right)  $, two-sided
p-value $p_{i}=2\Phi\left(  -\tilde{\sigma}_{ii}^{-1/2}\left\vert \zeta
_{i}\right\vert \right)  $, $X_{i}=1_{\left\{  p_{i}\leq t\right\}  }$,
$R_{m}\left(  t\right)  =\sum_{i=1}^{m}X_{i}$ and $V_{m}\left(  t\right)
=\sum_{i\in Q_{0,m}}X_{i}$.

We aim to show that the variance $\mathbb{V}\left[  m^{-1}R_{m}\left(
t\right)  \right]  $ of $m^{-1}R_{m}\left(  t\right)  $ satisfies $O\left(
m^{-\delta_{\ast}}\right)  $ with $\delta_{\ast} = \min \left\{\delta,1\right\}$. For a one-sided p-value $p_{i}$, define $\tilde
{t}=\tilde{\sigma}_{ii}^{1/2}\Phi^{-1}\left(  1-t\right)  $, $r_{1,i}=\tilde{t}%
-\mu_{i}$ and $r_{2,i}=-\infty$; for a two-sided p-value $p_{i}$, define
$\tilde{t}=-\tilde{\sigma}_{ii}^{1/2}\Phi^{-1}\left(  2^{-1}t\right)  $,
$r_{1,i}=\tilde{t}-\mu_{i}$ and $r_{2,i}=-\tilde{t}-\mu_{i}$. Further, set
$c_{l,i}=\tilde{\sigma}_{ii}^{-1/2}r_{l,i}$ for $l=1,2$, let $\rho_{ij}$ be the
correlation between $\zeta_{i}$ and $\zeta_{j}$ for $i\neq j$, and define the
sets%
\[
\left\{
\begin{array}
[c]{c}%
E_{1,m}=\left\{  \left(  i,j\right)  :1\leq i,j\leq m,i\neq j,\left\vert
\rho_{ij}\right\vert <1\right\}  ,\\
E_{2,m}=\left\{  \left(  i,j\right)  :1\leq i,j\leq m,i\neq j,\left\vert
\rho_{ij}\right\vert =1\right\}  .
\end{array}
\right.
\]
Namely, $E_{2,m}$ records pairs $\left(  \zeta_{i},\zeta_{j}\right)  $ with
$i\neq j$ such that $\zeta_{i}$ and $\zeta_{j}$ are linearly dependent. Obviously,
$\left\vert\mathrm{cov}\left(  X_{i},X_{j}\right) \right\vert \le C  = C \left\vert \rho_{ij}\right\vert$
for $\left(i,j\right) \in E_{2,m}$.
Further, %
\begin{align}
&  \mathbb{V}\left[  m^{-1}R_{m}\left(  t\right)  \right]  \nonumber\\
&  =m^{-2}\sum_{i=1}^{m}\mathbb{V}\left[  X_{i}\right]  +m^{-2}\sum_{\left(
i,j\right)  \in E_{2,m}}\mathrm{cov}\left(  X_{i},X_{j}\right)  +m^{-2}%
\sum_{\left(  i,j\right)  \in E_{1,m}}\mathrm{cov}\left(  X_{i},X_{j}\right)
.\label{eqE1}
\end{align}
However, we have $m^{-2}\sum_{i=1}^{m}\mathbb{V}\left[  X_{i}\right]
=O\left(  m^{-1}\right)  $ since the $\vert X_i \vert$'s as uniformly bounded by $1 $ a.s., and
\begin{align*}
m^{-2}\left\vert \sum\nolimits_{\left(  i,j\right)  \in E_{2,m}}%
\mathrm{cov}\left(  X_{i},X_{j}\right)  \right\vert  &  =O\left(
m^{-2}\left\vert E_{2,m}\right\vert \right)
  = O\left(  m^{-2}\left\Vert \mathbf{R}\right\Vert _{1}\right)  =O\left(
m^{-\delta}\right)  .
\end{align*}
This implies%
\begin{equation}
\left\vert m^{-2}\sum_{i=1}^{m}\mathbb{V}\left[  X_{i}\right]  +m^{-2}%
\sum\nolimits_{\left(  i,j\right)  \in E_{2,m}}\mathrm{cov}\left(  X_{i}%
,X_{j}\right)  \right\vert =O\left(  m^{-\min\left\{  \delta,1\right\}
}\right)  .\label{eqE11}%
\end{equation}
So, we only need to deal with $m^{-2}\sum_{\left(  i,j\right)  \in E_{1,m}%
}\mathrm{cov}\left(  X_{i},X_{j}\right)  $ in the right hand side (RHS) of
(\ref{eqE1}).

Consider first one-sided p-values and pick a pair $\left(  i,j\right)  \in
E_{1,m}$. By the definition of covariance,
\begin{align}
\mathrm{cov}\left(  X_{i},X_{j}\right)   &  =\int_{-\infty}^{c_{1,i}}%
\int_{-\infty}^{c_{1,j}}\left[  f_{\rho_{ij}}\left(  x,y\right)  -\phi\left(
x\right)  \phi\left(  y\right)  \right]  dxdy \nonumber \\
&  =\int_{-\infty}^{c_{1,i}}\int_{-\infty}^{c_{1,j}}\sum_{n=1}^{\infty}%
\frac{\rho_{ij}^{n}}{n!}H_{n}\left(  x\right)  H_{n}\left(  y\right)
dxdy,\label{eqE2B}%
\end{align}
where we have used Mehler expansion%
\begin{equation}
f_{\rho}\left(  x,y\right)  =\left(  1+\sum_{n=1}^{\infty}\frac{\rho^{n}}%
{n!}H_{n}\left(  x\right)  H_{n}\left(  y\right)  \right)  \phi\left(
x\right)  \phi\left(  y\right)  \label{eqE3}%
\end{equation}
for $\left\vert \rho_{ij}\right\vert \neq1$, $\phi\left(  x\right)
=\left(  2\pi\right)  ^{-1/2}\exp\left(  -x^{2}/2\right)  $, and the $n$th Hermite polynomial%

\[
H_{n}\left(  x\right)  =\left(  -1\right)  ^{n}\frac{1}{\phi\left(  x\right)
}\frac{d^{n}}{dx^{n}}\phi\left(  x\right)  .
\]
Since $\mathrm{cov}\left(  X_{i},X_{j}\right)  $ is well-defined, the RHS of
(\ref{eqE2B}) is convergent. However, by \cite{Watson:1933}, the series on the
RHS of (\ref{eqE3}) as a trivariate function of $\left(  x,y,\rho\right)  $ is
uniformly convergent on each compact set of $\mathbb{R}\times\mathbb{R}%
\times\left(  -1,1\right)  $. Therefore, we can interchange the order of
summation and integration on the RHS of (\ref{eqE2B}) to obtain%
\begin{equation}
\mathrm{cov}\left(  X_{i},X_{j}\right)  =\sum_{n=1}^{\infty}\int_{-\infty
}^{c_{1,i}}\int_{-\infty}^{c_{1,j}}\frac{\rho_{ij}^{n}}{n!}H_{n}\left(
x\right)  H_{n}\left(  y\right)  dxdy.\label{eqE4}%
\end{equation}
Plugging into (\ref{eqE4}) the identity
$
H_{n-1}\left(  x\right)  \phi\left(  x\right)  =\int_{-\infty}^{x}H_{n}\left(
y\right)  \phi\left(  y\right)  dy
$
for $x\in\mathbb{R}$ and $n\geq1$, we have%
\[
\mathrm{cov}\left(  X_{i},X_{j}\right)  =\sum_{n=1}^{\infty}\frac{\rho
_{ij}^{n}}{n!}H_{n-1}\left(  c_{1,i}\right)  H_{n-1}\left(  c_{1,j}\right)
\phi\left(  c_{1,i}\right)  \phi\left(  c_{1,j}\right)  .
\]

Now consider two-sided p-values and pick a pair $\left(  i,j\right)  \in
E_{1,m}$. Then%
\[
\mathrm{cov}\left(  X_{i},X_{j}\right)  =\int_{c_{2,i}}^{c_{1,i}}\int%
_{c_{2,j}}^{c_{1,j}}\sum_{n=1}^{\infty}\frac{\rho_{ij}^{n}}{n!}H_{n}\left(
x\right)  H_{n}\left(  y\right)  dxdy.
\]
Following the previous arguments for the case of one-sided p-values, we obtain%
\[
\mathrm{cov}\left(  X_{i},X_{j}\right)  =\sum_{n=1}^{\infty}\frac{\rho_{ij}^{n}}{n!}d_{n-1}\left(
c_{1,i},c_{2,i}\right)  d_{n-1}\left(  c_{1,j},c_{2,j}\right)  ,
\]
where $d_{n}\left(  c,c^{\prime}\right)  =H_{n}\left(  c\right)  \phi\left(
c\right)  -H_{n}\left(  c^{\prime}\right)  \phi\left(  c^{\prime}\right)$ for $c,c^{\prime}\in\mathbb{R}$.%

By \autoref{BndHermitePoly}, i.e.,%
\begin{equation}
\left\vert e^{-y^{2}/2}H_{n}\left(  y\right)  \right\vert \leq K_{0}\sqrt
{n!}n^{-1/12}e^{-y^{2}/4}\text{ \ for any\ }y\in\mathbb{R},\label{eqHermiteE}%
\end{equation}
and the uniform boundedness of $\phi$, we have%
\[
\left\vert \left(  n!\right)  ^{-1}H_{n-1}\left(  c_{1,i}\right)
H_{n-1}\left(  c_{1,j}\right)  \phi\left(  c_{1,i}\right)  \phi\left(
c_{1,j}\right)  \right\vert \leq Cn^{-7/6}%
\]
and%
\[
\left\vert \left(  n!\right)  ^{-1}d_{n-1}\left(  c_{1,i},c_{2,i}\right)
d_{n-1}\left(  c_{1,j},c_{2,j}\right)  \right\vert \leq Cn^{-7/6}%
\]
for a generic constant $C>0$. So, for each $\left(  i,j\right)  \in E_{1,m}$ and both
types of p-values,%
\begin{equation*}
\left\vert \mathrm{cov}\left(  X_{i},X_{j}\right)  \right\vert \leq C\left\vert
\rho_{ij}\right\vert \sum_{n=1}^{\infty}n^{-7/6}%
\end{equation*}
and%
\begin{align}
\left\vert m^{-2}\sum_{\left(  i,j\right)  \in E_{1,m}}\mathrm{cov}\left(
X_{i},X_{j}\right)  \right\vert  & \leq Cm^{-2}\sum_{\left(  i,j\right)  \in
E_{1,m}}\left\vert \rho_{ij}\right\vert \sum_{n=1}^{\infty}n^{-7/6}%
\nonumber\\
& \leq Cm^{-2}\sum_{\left(  i,j\right)  \in E_{1,m}}\left\vert \rho
_{ij}\right\vert \le Cm^{-2}\left\Vert \mathbf{R}\right\Vert _{1} =O\left(  m^{-\delta}\right) \label{eqE13}
\end{align}
by also observing that $\sum_{n=1}^{\infty}n^{-7/6}$ is convergent. Combining
(\ref{eqE1}), (\ref{eqE11}) and (\ref{eqE13}), we have%
\begin{equation}
\mathbb{V}\left[  m^{-1}R_{m}\left(  t\right)  \right]  =O\left(
m^{-\min\left\{  \delta,1\right\}  }\right)    \label{eqVarPBE}%
\end{equation}
for both types of p-values. Identical
arguments assert that $\mathbb{V}\left[  m^{-1}V_{m}\left(  t\right)  \right]
=O\left(m^{-\min\left\{  \delta,1\right\}  }\right)  $. Thus, the claims on
$ m^{-1}R_{m}\left(  t\right)  $ and $m^{-1}V_{m}\left(  t\right) $ are valid.
Further, when $\liminf_{m \to \infty}\pi_{0,m}>0$,
\begin{equation}\label{SLLNVm}
  \left \vert m_0^{-1} V_{m}\left(t\right) -\mathbb{E}\left[  m_0^{-1}V_{m}\left(  t\right) \right] \right \vert \to 0 \quad \text{a.s.}
\end{equation}

Secondly, we show the second claim. Let us assume for the moment
\begin{equation}\label{lastCond}
  \liminf_{m\rightarrow \infty}m^{-1}R_{m}\left(  t\right)  >0 \quad \text{a.s.}
\end{equation}
 and prove \eqref{lastCond} at the end. Setting $r_{\ast}=\liminf
_{m\rightarrow\infty}m^{-1}R_{m}\left(  t\right)  $. Then $r_{\ast}>0$ a.s.,
\begin{equation}
\liminf_{m\rightarrow\infty}\mathbb{E}\left[  m^{-1}R_{m}\left(  t\right)
\right]  \geq r_{\ast}>0,\label{eqE6}%
\end{equation}
and
\begin{equation}
R_{m}\left(  t\right)  \vee1=R_{m}\left(  t\right)  \text{ for all }m\text{
large enough}.\label{eqE15}%
\end{equation}
Recall $\mathrm{FDP}_{m}\left(  t\right)  = \frac{V_{m}\left(  t\right)  }%
{R_{m}\left(  t\right) \vee 1}$. So, \eqref{eqE15} implies
\begin{equation*}
\mathrm{FDP}_{m}\left(  t\right)  =\frac{m^{-1}V_{m}\left(  t\right)  }{m^{-1}%
R_{m}\left(  t\right)  }
\end{equation*}
for all $m$ large enough, and with \eqref{eqE6} the continuous mapping theorem implies a.s.%
\begin{equation}
\lim_{m\rightarrow\infty}\left\vert \mathrm{FDP}_{m}\left(  t\right)  -\frac
{\mathbb{E}\left[  m^{-1}V_{m}\left(  t\right)  \right]  }{\mathbb{E}\left[
m^{-1}R_{m}\left(  t\right)  \right]  }\right\vert =0.\label{eqE8}%
\end{equation}
Since by \eqref{eqE6} a.s.%
\begin{equation}
\left\vert \mathrm{FDP}_{m}\left(  t\right)  -\frac{\mathbb{E}\left[  m^{-1}%
V_{m}\left(  t\right)  \right]  }{\mathbb{E}\left[  m^{-1} \left(R_{m}\left(
t\right)  \vee1\right)\right]  }\right\vert \leq2,\label{eqE9}%
\end{equation}
applying the dominated convergence theorem together with \eqref{eqE15}, \eqref{eqE8} and (\ref{eqE9}) gives%
\begin{equation}
\lim_{m\rightarrow\infty}\mathbb{E}\left[  \mathrm{FDP}_{m}\left(  t\right)
-\frac{\mathbb{E}\left[  m^{-1}V_{m}\left(  t\right)  \right]  }%
{\mathbb{E}\left[  m^{-1}R_{m}\left(  t\right)   \right]  }\right]
=0.\label{eqE10}%
\end{equation}
Applying \eqref{eqE15}, (\ref{eqE8}) and (\ref{eqE10}) to the decomposition%
\begin{align*}
& \mathrm{FDP}_{m}\left(  t\right)  -\mathbb{E}\left[  \mathrm{FDP}_{m}\left(  t\right)
\right]  \\
& =\mathrm{FDP}_{m}\left(  t\right)  -\frac{\mathbb{E}\left[  m^{-1}V_{m}\left(
t\right)  \right]  }{\mathbb{E}\left[  m^{-1}R_{m}\left(  t\right)
 \right]  }-\left(  \mathbb{E}\left[  \mathrm{FDP}_{m}\left(  t\right)  \right]
-\frac{\mathbb{E}\left[  m^{-1}V_{m}\left(  t\right)  \right]  }%
{\mathbb{E}\left[  m^{-1}R_{m}\left(  t\right)  \right]  }\right)  \\
& =\mathrm{FDP}_{m}\left(  t\right)  -\frac{\mathbb{E}\left[  m^{-1}V_{m}\left(
t\right)  \right]  }{\mathbb{E}\left[  m^{-1}R_{m}\left(  t\right)
 \right]  }-\mathbb{E}\left[  \mathrm{FDP}_{m}\left(  t\right)  -\frac
{\mathbb{E}\left[  m^{-1}V_{m}\left(  t\right)  \right]  }{\mathbb{E}\left[
m^{-1}R_{m}\left(  t\right)   \right]  }\right]
\end{align*}
gives a.s.%
\[
\lim_{m\rightarrow\infty}\left\vert \mathrm{FDP}_{m}\left(  t\right)  -\mathbb{E}%
\left[  \mathrm{FDP}_{m}\left(  t\right)  \right]  \right\vert =0.
\]

Finally, we show \eqref{lastCond}, i.e., $\liminf_{m\rightarrow \infty}m^{-1}R_{m}\left(  t\right)  >0$ a.s. when $\liminf_{m \to \infty}\pi_{0,m}>0$. Since the p-values associated with the true null hypotheses $\mu_i =0$ for $i \in Q_{0,m}$ are identically distributed, $t \in \left(0,1\right)$ is a fixed, positive constant, and \eqref{SLLNVm} holds, we have
\begin{align*}
\liminf_{m\rightarrow\infty}\frac{R_{m}\left(  t\right)  }{m} &  \geq
\liminf_{m\rightarrow\infty}\frac{m_{0}}{m}\frac{V_{m}\left(  t\right)
}{m_{0}}\geq\liminf_{m\rightarrow\infty} \pi_{0,m} \times
\liminf_{m\rightarrow\infty}\frac{V_{m}\left(  t\right)  }{m_{0}}\\
&  \geq C\frac{m_{0}\mathbb{E}\left[  1_{\left\{  p_{0}\leq t\right\}
}\right]  }{m_{0}}=C\mathbb{E}\left[  1_{\left\{  p_{0}\leq
t\right\}  }\right]  > 0 \quad \text{a.s.}
\end{align*}
for some constant $C>0$, where $p_{0}$ is the p-value associated with a true null hypothesis.
This completes the proof.

\subsection{Proof of \autoref{ThmSLLNAMarg}}

Note that $\ms{\eta} = \mathbf{0}$ a.s. So, the set $G_{m,\boldsymbol{\eta}}\left(  t,\varepsilon_{m}\right)$ defined in \eqref{eqExceptionB} changes into
\begin{equation*}
G_{m}\left(  t,\varepsilon_{m}\right)  =\bigcup
\nolimits_{i\in E_{0,m}}\left\{\mu_i,\tilde{\sigma}_{ii}:\min\left\{
\left\vert r_{1,i}\right\vert ,\left\vert r_{2,i}\right\vert \right\}
<\varepsilon_{m}\right\},
\end{equation*}
where $\tilde{\sigma}_{ii}$ is the standard deviation of $\zeta_i$ and $r_{ij}$'s are defined in the beginning of \autoref{secVarianceOfCondRej}.
Further, from \autoref{ThmSLLN}, we see the following: if $\sigma_{0}>0$, then%
\begin{equation*}
\mathbb{V}\left[  m^{-1}R_{m}\left(  t  \right)  \right]  =O\left(  m^{-\min\left\{
\delta,1\right\}  }\right)  ; 
\end{equation*}
otherwise, for each pair $\rbk{\ms{\mu},\ms{\sigma}}\notin G_{m}\left(  t,\varepsilon_{m}\right)$,
\begin{equation*}
\mathbb{V}\left[  m^{-1}R_{m}\left(  t  \right)  \right] =O\left(  \varepsilon
_{m}^{-2}m^{-\min\left\{  \delta,1\right\}  }\right) =  O\left(  m^{-\delta^{\prime} }\right) \ \text{for some \ } \delta^{\prime}>0.
\end{equation*}

Hence, by \autoref{LemmaLyons}, the
conclusions hold for $R_{m}$. On the other hand, as can be seen from the proof of \autoref{ThmSLLN}, the upper bounds on the variance of $m^{-1}R_{m}\left(  t\right)$ are also upper bounds on the variance of $m^{-1}V_{m}\left(  t\right)$ for each $m \ge 1$. So, the assertions on $V_{m}\left(  t\right)$ are valid. Further, $m_0^{-1}V_{m}\left(  t\right) = m m_0^{-1}  m^{-1}V_{m}\left(  t\right)$. So,
$\left \vert m_0^{-1} V_{m}\left(t\right) -\mathbb{E}\left[  m_0^{-1}V_{m}\left(  t\right) \right] \right \vert \to 0$ a.s.
when $\liminf_{m \to \infty}\pi_{0,m}>0$. Finally, the assertions on $\left\vert \mathrm{FDP}_{m}\left(  t\right)  -\mathbb{E}\left[
\mathrm{FDP}_{m}\left(  t  \right)  \right]  \right\vert \to 0$ a.s. can be proved using the same arguments in the second and third parts of the proof of \autoref{prop:SLLNMarginal}. This completes the proof.

\subsection{Proof of \autoref{ThmSLLNStrongest}}
\label{secPfsllna}

Whenever needed, we will also write $v_{i}$, $\eta_{i}$, $\mu_{i}$, $\zeta
_{i}$ and $X_{i}$ with $1\leq i\leq m$ as $v_{i,m}$, $\eta_{i,m}$ $\mu_{i,m}$,
$\zeta_{i,m}$ and $X_{i,m}$ for $m\geq1$. The proof is divided into
\textbf{Step 1} to show the SLLN for a subsequence of $\left\{  X_{i}\right\}
$ and \textbf{Step 2} to show the SLLN for the sequence $\left\{
X_{i}\right\}  $ via a controlled maximal inequality.

\textbf{Step 1}: Recall $X_{i,m}=1_{\left\{  p_{i}\leq t|\eta_{i,m}\right\}
}$. By \autoref{ThmSLLN}, we have
\[
\left\vert \mathrm{cov}_{\mathbf{v}}\left(  X_{i,m},X_{j,m}\right)
\right\vert \leq C\left\vert \rho_{ij}\right\vert \text{ \ for all \ }1\leq
i\leq j\leq m,
\]
where $\rho_{ij}$ is the correlation between $v_{i,m}$ and $v_{j,m}$. Let
$\tilde{X}_{i,m}=X_{i,m}-\mathbb{E}_{\mb{v}}\left[  X_{i,m}\right]  $ and
$\mathbf{K}_{m}$ be the covariance matrix of $\left\{  \tilde{X}_{i}\right\}
_{i=1}^{m}$. Then, $m^{-2}\left\Vert \mathbf{R}_{\mathbf{v}}\right\Vert
_{1}=O\left(  m^{-\delta}\right)  $ implies%
\[
\mbb{V}_{\mb{v}}\left(  m^{-1}\sum\nolimits_{i=1}^{m}\tilde{X}_{i,m}\right)  \leq
Cm^{-2}\left\Vert \mathbf{K}_{m}\right\Vert _{1}\leq Cm^{-2}\left\Vert
\mathbf{R}_{\mathbf{v}}\right\Vert _{1}=O\left(  m^{-\delta}\right)  .
\]
Let $b_{m}=m^{-2}\left\Vert \mathbf{R}_{\mathbf{v}}\right\Vert _{1}$. Then
$b_{m}=O\left(  m^{-\delta}\right)  $ and $\sum_{m\geq1}m^{-1}b_{m}<\infty$.
So, the lemma in \cite{Dvoretzky:1949} (which is restated as Lemma 2 of
\cite{Lyons:1988}) implies $\sum_{k\geq1}b_{m_{k}}<\infty$ for a subsequence
$m_{k}$ such that $m_{k}\rightarrow\infty$ and $m_{k+1}/m_{k}\rightarrow1$.
Therefore, Lemma 3 of \cite{Lyons:1988} implies that $m_{k}^{-1}\sum
_{i=1}^{m_{k}}\tilde{X}_{i,m_{k}}\rightarrow0$ a.s.

\textbf{Step 2}: Let $m$ be such that $m_{k}\leq m<m_{k+1}$. Then%
\begin{align*}
\left\vert \frac{1}{m}\sum\nolimits_{i=1}^{m}\tilde{X}_{i,m}\right\vert  &
\leq\frac{1}{m_{k}}\left\vert \sum\nolimits_{i=1}^{m_{k}}\left(  \tilde
{X}_{i,m}-\tilde{X}_{i,m_{k}}\right)  \right\vert +\left\vert \frac{1}{m_{k}%
}\sum\nolimits_{i=1}^{m_{k}}\tilde{X}_{i,m_{k}}\right\vert \\
&  +\max_{1\leq s\leq m-m_{k}}\left\vert \frac{1}{m_{k}}\sum\nolimits_{i=m_{k}%
+1}^{m_{k}+s}\tilde{X}_{i,m}\right\vert ,
\end{align*}
for which on the right-hand side the third term converges to $0$ a.s. since
$\left\vert \tilde{X}_{i,m}\right\vert \leq1$ a.s. and $m_{k+1}%
/m_{k}\rightarrow1$, and the second term converges to $0$ a.s. as already
justified. So, it is left to show%
\begin{equation}
\frac{1}{m_{k}}\left\vert \sum\nolimits_{i=1}^{m_{k}}\left(  \tilde{X}%
_{i,m}-\tilde{X}_{i,m_{k}}\right)  \right\vert \rightarrow0\text{ a.s.}
\label{eqab1}%
\end{equation}
Recall the notations: for a one-sided p-value $p_{i}$, define $\tilde{t}=\tilde{\sigma
}_{ii}^{1/2}\Phi^{-1}\left(  1-t\right)  $, $r_{1,i}=\tilde{t}-\mu_{i}%
-\eta_{i}$ and $r_{2,i}=-\infty$; for a two-sided p-value $p_{i}$, define
$\tilde{t}=-\tilde{\sigma}_{ii}^{1/2}\Phi^{-1}\left(  2^{-1}t\right)  $,
$r_{1,i}=\tilde{t}-\mu_{i}-\eta_{i}$ and $r_{2,i}=-\tilde{t}-\mu_{i}-\eta_{i}%
$; set $c_{l,i}=\sigma_{i,m}^{-1}r_{l,i}$ for $l=1,2$. Write $r_{l,i}$ and
$c_{l,i}$ respectively as $r_{l,i,m}$ and $c_{l,i,m}$. Let $B_{i,m}  =\left\{  \omega\in\Omega:X_{i,m}=0\right\}  $. Then%
\[
B_{i,m}  =\left\{  \omega\in\Omega:c_{2,i,m}%
\leq\sigma_{i,m}^{-1}v_{i,m}\leq c_{1,i,m}\right\}
\]
and $d_{i,m}=\mathbb{P}\left(  X_{i,m}=0\right)  =\Phi\left(  c_{1,i,m}\right)
-\Phi\left(  c_{2,i,m}\right)  $. Recall $\sigma_{0}=\lim_{m\rightarrow\infty
}\min_{1\leq i\leq m}\left\{  \sigma_{i,m}:\sigma_{i,m}\neq0\right\}  $ and pick any $m,m^{\prime}$ such that
$m\leq m^{\prime}$ and $m\rightarrow\infty$. However, $\sigma
_{0}>0$, $\max_{1\leq i\leq m}\left\vert \mu_{i,m}-\mu_{i,m^{\prime}}\right\vert
\rightarrow0$, $\max_{1\leq i\leq m}\left\vert \eta_{i,m}-\eta_{i,m^{\prime}%
}\right\vert \rightarrow0$ a.s. and $\max_{1\leq i\leq m}\left\vert
\sigma_{i,m}-\sigma_{i,m^{\prime}}\right\vert \rightarrow0$, So,
$B_{i,m}  \ominus B_{i,m^{\prime}}  \rightarrow\varnothing$ a.s. uniformly in $1\leq
i\leq m$ where $\ominus$ is the symmetric set difference, and $\left\vert
d_{i,m}  -d_{i,m^{\prime}}  \right\vert \rightarrow0$ uniformly in $1\leq i\leq m$ where we
have used the continuity of $\Phi$. Therefore, $$\max_{1\leq i\leq m_{k}%
}\left\vert \tilde{X}_{i,m}-\tilde{X}_{i,m_{k}}\right\vert \rightarrow0 \quad \text{a.s.}$$
and (\ref{eqab1}) holds. Consequently,$\ m^{-1}\sum\nolimits_{i=1}^{m}%
\tilde{X}_{i,m}\rightarrow0$ a.s.

Following the arguments given above, we can prove that $m^{-1}\sum
\nolimits_{i\in Q_{0,m}}\tilde{X}_{i,m}\rightarrow0$ a.s., where $Q_{0,m}$ contains $i$ such that
$1 \le i \le m$ and $\mu_{i}=0$. Since $V_{m}\left(  t|\boldsymbol{\eta}\right)
=\sum\nolimits_{i\in Q_{0,m}}{X}_{i,m}$, the assertion on $V_{m}\left(  t|\boldsymbol{\eta}\right)$ holds. Following the arguments in the last two parts
of the proof of \autoref{prop:SLLNMarginal}, we can show $$m^{-1}\left\vert \mathrm{FDP}%
_{m}\left(  t|\boldsymbol{\eta}\right)  -\mathbb{E}_{\mathbf{v}}\left[
\mathrm{FDP}_{m}\left(  t|\boldsymbol{\eta}\right)  \right]  \right\vert
\rightarrow0 \quad \text{a.s.}$$ when in addition $\liminf_{m\rightarrow\infty}\pi_{0,m}>0$.
However, we omit the remaining details here. This completes the proof.

\section{SLLN associated with adjusted conditional MTP}
\label{secAcmtp}

For the adjusted conditional MTP based on PCS, rather than providing complete proofs we will just point out the differences in the conditions and arguments, if any, that lead to the corresponding SLLN. In this section, we will maintain the same notations used in
\autoref{secNormalMeans} and point out the differences in their meanings if any.
Now abbreviate ``adjusted conditional MTP'' as ``acMTP''. For the acMTP, the one-sided p-value is
$\tilde{p}_i = 1-\Phi\big(\sigma_{i,m}^{-1}\left( \zeta_i - \eta_i \right)\big)$, and the two-sided
$\tilde{p}_i = 2 \Phi\big(-\sigma_{i,m}^{-1}\vert \zeta_i - \eta_i \vert \big)$
by observing $\zeta_i - \eta_i = \mu_i + v_i$ and $\mu_i + v_i \sim \mathsf{N}_{1}\big( \mu_i,\sigma_{i,m}^{2}\big)$.
When $\ms{\eta} = \mb{0}$ a.s., the acMTP is just the marginal MTP.

For the acMTP, let $X_{i}=1_{\left\{  \tilde{p}_{i}\leq t\right\}
}$ be the indicator of whether $\tilde{p}_{i}$ is
no larger than $t$. Then $X_{i}= 1_{\left\{  \tilde{p}_{i}\leq t |\boldsymbol{\eta}\right\}}$,
and the acMTP rejects $H_{i0}: \mu_i = 0$ iff $\tilde{p}_i \le t$.
Set $R_{m}\left(t|\boldsymbol{\eta}\right)  =\sum_{i=1}^{m}X_{i}$. The key to derive the SLLN
for $\left\{  R_{m}\left(  t|\boldsymbol{\eta}\right)  \right\}  _{m\geq1}$ is
to obtain the variance $\mathbb{V}_{\mathbf{v}}\left[  m^{-1}R_{m}\left(
t|\boldsymbol{\eta}\right)  \right]  $ for $m^{-1}R_{m}\left(
t|\boldsymbol{\eta}\right)  $ by expanding $\mathbb{V}_{\mathbf{v}}\left[
m^{-1}R_{m}\left(  t|\boldsymbol{\eta}\right)  \right]  $ into summands each
being an integral. Compared to the derivation for the variance the average
number of rejections for the conditional MTP in \autoref{secNormalMeans}, the only difference
are the changes in the upper and lower limits in the summands that made up
$\mathbb{V}_{\mathbf{v}}\left[  m^{-1}R_{m}\left(  t|\boldsymbol{\eta}\right)
\right]  $. Specifically, we only have to change the $\tilde{t}$'s defined in
the beginning of \autoref{secVarianceOfCondRej} into $\tilde{t}=\sigma_{i,m}\Phi^{-1}\left(  1-t\right)  $ for a
one-sided p-value $\tilde{p}_{i}$ or $\tilde{t}=-\sigma_{i,m}\Phi^{-1}\left(
2^{-1}t\right)  $ for a two-sided p-value $\tilde{p}_{i}$, and maintain the
definitions of all other quantities (with $X_i$ taking the new meaning here).
So, the arguments that lead to \autoref{ThmSLLNStrongest} lead
to the following:
\begin{proposition}\label{SLLNacMTP}
Consider the acMTP. Assume $\sigma_0>0$ and
   $m^{-2}\left\Vert \mathbf{R}_{\mathbf{v}}\right\Vert _{1} = O\left(   m^{-\delta
} \right)$ for some $\delta >0$.
If for any natural numbers
$m,m^{\prime}$ such that $m\leq m^{\prime}$ and $m\rightarrow\infty$%
\begin{equation*}
\max\left\{  \left\Vert \boldsymbol{\eta}_{m}-\boldsymbol{\eta}_{m^{\prime}%
}^{\left(  m\right)  }\right\Vert _{2},\left\Vert \boldsymbol{\mu}%
_{m}-\boldsymbol{\mu}_{m^{\prime}}^{\left(  m\right)  }\right\Vert
_{2},\left\Vert \boldsymbol{\sigma}_{\mathbf{v},m}-\boldsymbol{\sigma
}_{\mathbf{v},m^{\prime}}^{\left(  m\right)  }\right\Vert _{2}\right\}
\rightarrow0\text{ a.s.,}%
\end{equation*}
then
$m^{-1}\left\vert R_{m}\left(  t|\ms{\eta}  \right)  -\mathbb{E}_{\mathbf{v}}\left[  R_{m}\left(
t|\ms{\eta}  \right)  \right]
\right\vert \to 0$ a.s and $m^{-1}\left\vert V_{m}\left(  t|\ms{\eta}  \right)  -\mathbb{E}_{\mathbf{v}}\left[  V_{m}\left(
t|\ms{\eta}  \right)  \right] \right\vert \to 0$ a.s.
If in addition $\liminf_{m\rightarrow\infty}\pi_{0,m}  >0$ and $\lim_{m \to \infty}\min_{i \in Q_{0,m}} \tilde{p}_i >0$ a.s. uniformly in $\ms{\eta}$,  then
$m^{-1}\left\vert \mathrm{FDP}_{m}\left(  t|\ms{\eta}\right)
 -\mathbb{E}_{\mathbf{v}}\left[  \mathrm{FDP}_{m}\left(  t|\ms{\eta}\right)  \right]
\right\vert \to 0$ a.s.
\end{proposition}

In the statement above, $\mathbf{R}_{\mathbf{v}}$ is the correlation matrix of the minor vector $\mathbf{v}$. We remark that the extra condition ``$\lim_{m \to \infty}\min_{i \in Q_{0,m}} \tilde{p}_i >0$ a.s. uniformly in $\ms{\eta}$'' is used to ensure $\liminf_{m \to \infty} m^{-1}R_{m}\left(  t|\ms{\eta}  \right) >0$ a.s., so that the continuous mapping theorem can be applied to show the assertion on  $m^{-1}\mathrm{FDP}_{m}\left(  t|\ms{\eta}\right)$.

\section{Relationship between PCS and PFA}
\label{SecConnection}

In general, PCS is different than PFA, as we argue as follows.
Let $\left\{  \lambda_{i,m}\right\}  _{i=1}^{m}$ be the
descendingly ordered (in $i$) eigenvalues (counting multiplicity) of
$\boldsymbol{\Sigma}$ whose corresponding eigenvectors are $\boldsymbol{\gamma
}_{i}=\left(  \gamma_{i1},...,\gamma_{im}\right)  ^{\top}$ for $1\leq i\leq m$.
For some integer $k$ between $1$ and $m$, setting $\mathbf{w}=\left(
w_{1},...,w_{m}\right)  ^{\top}\sim\mathsf{N}_{m}\left(  \mathbf{0}%
,\mathbf{I}\right)  $ with $\mathbf{I}$ being the identity matrix,
\[
\boldsymbol{\eta}=\sum_{j=1}^{k}\lambda_{j,m}^{1/2}\boldsymbol{\gamma}%
_{j}w_{j}\text{ \ \ \ and \ \ \ }\mathbf{v}=\sum_{j=k+1}^{m}\lambda
_{j,m}^{1/2}\boldsymbol{\gamma}_{j}w_{j}%
\]
gives (\ref{eqModel}), i.e., $\boldsymbol{\zeta}=\boldsymbol{\mu
}+\boldsymbol{\eta}+\mathbf{v}$.
Recall $\boldsymbol{\Sigma}_{\mathbf{v}}=\left(  q_{ij}\right)  _{m\times m}$
is the covariance matrix of $\mathbf{v}$.
When $\ms{\Sigma}$ is a correlation matrix, set $\vartheta_{m}=m^{-1}\sqrt{\sum_{i=k+1}^{m}\lambda_{i,m}^{2}}$.
Then $\vartheta_{m}=m^{-1}\left\Vert \boldsymbol{\Sigma}_{\mathbf{v}}\right\Vert
_{2}$. Pick a
$\delta\in(0,1]$ and assume the existence of the smallest $k=k\left(
\delta,m\right)  $ between $1$ and $m$ such that $\vartheta_{m}=O\left(
m^{-\delta}\right)  $. Then the corresponding decomposition $\boldsymbol{\zeta
}=\boldsymbol{\mu}+\boldsymbol{\eta}+\mathbf{v}$ is the PFA in \cite{Fan:2012},
where for simplicity we also refer to $\ms{\eta}$ and $\mb{v}$ as the major and minor vectors, respectively. In this case,
\begin{equation}
m^{-1}\left\Vert \boldsymbol{\Sigma}_{\mathbf{v}}\right\Vert
_{2} = O\left(m^{-\delta}\right), \label{eqFanCond}%
\end{equation}
which is the only condition used by \cite{Fan:2012} to claim the SLLN associated with the conditional MTP.
Since the inequality
$m^{-2}\left\Vert \boldsymbol{\Sigma}_{\mathbf{v}}\right\Vert _{1}\leq
m^{-1}\left\Vert \boldsymbol{\Sigma}_{\mathbf{v}}\right\Vert _{2}$
implies $m^{-2}\left\Vert \boldsymbol{\Sigma}_{\mathbf{v}}\right\Vert
_{1}=O\left(  m^{-\delta}\right)  $ but not $m^{-2}\left\Vert \mathbf{R}_{\mathbf{v}}\right\Vert
_{1}=O\left(  m^{-\delta}\right)  $ (as the definition of PCS), in general PCS is different from PFA.
On the other hand, it is easily seen that PCS can be realized by PFA when the variances of $\mb{v}$ are uniformly bounded away from $0$ and $\infty$ for all $m$.

\section*{Acknowledgements}

This research was funded in part by a National Science Foundation Plant Genome
Research Program grant (No. IOS-1025976) to R.W. Doerge. Part of it was
completed when X. Chen was a PhD candidate at Purdue University. We thank John D. Storey for very kind support,
William B. Johnson for pointing out Walsh matrices, and Jo\~{a}o M. Pereira for a discussion on constructing
orthogonal matrices with special properties using discrete Fourier transform.

\bibliographystyle{chicago}

\end{document}